\documentclass{article}
\pdfoutput=1
\usepackage[utf8]{inputenc}
\usepackage[english]{babel}
\usepackage{amsmath, amsthm, amssymb, amsfonts}
\usepackage{graphicx, color, placeins}

\usepackage{tikz-cd}

\usepackage{enumitem}
\usepackage{hyperref}

\usepackage{array}

\newcommand{\cc}[1]{\mathcal{#1}}
\newcommand{\bb}[1]{\mathbb{#1}}

 \def\cB{{\mathcal{B}}}                      \def\cX{{\mathcal{X}}}  

    \def\bE{{\mathbb{E}}}         \def\bN{{\mathbb{N}}}   \def\bQ{{\mathbb{Q}}} \def\bR{{\mathbb{R}}}        \def\bZ{{\mathbb{Z}}}

\theoremstyle{plain}
\newtheorem{thm}{Theorem}[section]
\newtheorem{lem}[thm]{Lemma}

\theoremstyle{definition}

\theoremstyle{remark}

\renewcommand{\leq}{\leqslant} 
\renewcommand{\geq}{\geqslant}
\renewcommand{\epsilon}{\varepsilon} 
\renewcommand{\subset}{\subseteq} 

\renewcommand{\{}{\lbrace}
\renewcommand{\}}{\rbrace}
\newcommand{\sm}{\setminus} 

\renewcommand{\bar}{\overline}

\renewcommand{\P}{\cc{P}}
\newcommand{\Q}{\mathcal{Q}}
\newcommand{\E}{\mathbb{E}^3}
\newcommand{\TT}{\mathbb{T}^3}
\newcommand{\T}{\mathbf{T}}
\newcommand{\I}{\mathbf{I}}
\newcommand{\IE}{\I(\E)}
\renewcommand{\O}{\mathbf{O}}
\newcommand{\norm}{N_{\IE}(\bfL)}

\renewcommand{\L}{\Lambda}
\newcommand{\bfL}{\mathbf{\Lambda}}
\newcommand{\TTL}{\mathbb{T}^{3}_{\bfL}}

\renewcommand{\r}{\rho}

\newcommand{\+}{\oplus}

\newcommand{\tet}{\mathbf{[3,3]^+}}
\newcommand{\Tet}{\mathbf{[3,3]}}
\newcommand{\Oct}{\mathbf{[3,4]}}
\newcommand{\Ico}{\mathbf{[3,5]}}
\newcommand{\Coxpq}{\mathbf{[p,q]}}
\newcommand{\coxpq}{\mathbf{[p,q]^{+}}}
\newcommand{\DD}{\mathbf{\tilde{D}}}
\newcommand{\DDn}{\mathbf{\tilde{D}_{n}}}

\newcommand{\cl}{\L_{(1,0,0)}}
\newcommand{\fcl}{\L_{(1,1,0)}}
\newcommand{\bcl}{\L_{(1,1,1)}}
\newcommand{\tl}{\L^{\{3,6\}}_{(1,0)}}
\newcommand{\tcl}{\L^{\{3,6\}}_{(1,1)}}
\newcommand{\sql}{\L_{(1,0)}}
\newcommand{\scl}{\L_{(1,1)}}

\newcommand{\Cl}{\mathbf{\cl}}
\newcommand{\Fcl}{\mathbf{\fcl}}
\newcommand{\Bcl}{\mathbf{\bcl}}
\newcommand{\Tl}{\mathbf{\tl}}
\newcommand{\Tcl}{\mathbf{\tcl}}
\newcommand{\Sql}{\mathbf{\sql}}
\newcommand{\Scl}{\mathbf{\scl}}

\newcommand{\F}{\cc{F}}

\newcommand{\PL}{\P_{\bfL}}
\newcommand{\piL}{\pi_{\bfL}: \E \to \TTL}
\newcommand{\gop}{G_{o}(\P)}
\newcommand{\bgop}{\bar{G}_{o}(\P)}
\newcommand{\TB}{\cB_{\{3,6\}}}

\DeclareMathOperator{\diag}{Diag}
\DeclareMathOperator{\G}{\Gamma}
\DeclareMathOperator{\sk}{Sk^{1}}
\DeclareMathOperator{\bd}{bd}
\DeclareMathOperator{\clos}{cl}
\DeclareMathOperator{\inte}{int}
\DeclareMathOperator{\aff}{Aff}

\begin{document}

\title{Regular polyhedra in the $3$-torus.}

\author{Antonio Montero \thanks{Supported by CONACyT grant 414098 and PAPIIT UNAM project IN101615.}\\ amontero@matmor.unam.mx\\  Centro de Ciencias Matemáticas UNAM \\ Morelia, México. }
\date{\today}
\maketitle

\begin{abstract}
In this paper we discuss the classification rank $3$ lattices preserved by finite orthogonal groups of isometries and derive from it the classification of regular polyhedra in the $3$-dimensional torus. This classification is highly related to the classification of regular polyhedra in the $3$-space.
\end{abstract}

\section{Introduction} \label{sec:intro}
Symmetric structures such as polyhedra and tessellations had been of interest for centuries. Greeks were aware of platonic solids more that two thousand years ago and several similar structures appeared in art in the middle ages.

Regular polyhedra (those with maximum symmetry) have been one point where geometrical ideas fit together with combinatorics and group theory; an example of this is shown widely in McMullen and Schulte \cite{ARP}. However, the idea of treating polyhedra as both geometrical and combinatorial objects goes back before \cite{ARP}, being Coxeter one of the most remarkable exponents in XX century. Coxeter's work in regular polyhedra cannot be summarized in a few lines but it goes from work in his youth such as \cite{coxeterSkew} to its popular book \cite{coxeterRegularPolytopes} .

The idea of regular polyhedra has been modified several times and this has given birth to structures as those studied by Coxeter in \cite{coxeterSkew}. In 1977 Branko Grünbaum \cite{grunbaumOldAndNew} considered one of the most general definitions of regular polyhedra in Eclidean space $\E$ and gave a list of 47 structures. Grünbaum set aside the notions of finiteness, convexity or planarity of the faces (or vertex figures) and used a more graph-theoretical approach. A couple of years after later Andreas Dress completed the list to 48 regular polyhedra in Euclidean $3$-space and proved that the list was complete (see \cite{dress1,dress2}).

In 1997 Peter McMullen and Egon Schulte proved again that the list of Grünbaum-Dress regular polyhedra is complete starting with the concept of abstract polyhedron (see \cite{ordinary}). After the work of McMullen and Schulte there have been several generalizations. Many of them relax the symmetry condition (see \cite{pellWeissChiralR3,schulteChiralR31,schulteChiralR32} for instance) and higher ranks (dimensions) have been explored (see McMullen
\cite{mcmullenFourDimRegPolyh,mcmullenRegApeirotopesOfDimRank4,mcmullenRegularFullRank})

When changing the ambient space there is no much work published so far. In \cite{roliProjective1,roliProjective2} Bracho et al classified regular polyhedra with planar faces in projective space. The complete classification may be obtained as a consequence of McMullen's work in \cite{mcmullenFourDimRegPolyh}. Regular polyhedra (maps) in the $2$-torus are well known (see \cite[Chapter 8]{coxeterMoserGenandRelforDG}).

In Section\nobreakspace \ref {sec:basics} we review the basic notions of abstract polyhedra and regular polyhedra in $3$-space. In order to study isometries of the $3$-torus, in Section\nobreakspace \ref {sec:lattices} we classify rank 3 lattices preserved by finite orthogonal groups and finally in Section\nobreakspace \ref {sec:RPT} we give a complete classification of regular polyhedra in the $3$-torus.

\section{Basic Notions}	\label{sec:basics}

\subsection{Abstract Polyhedra}
 
Here we will introduce the notions of abstract polyhedra and their geometric realizations in Euclidean space. These are combinatorial generalizations of classical convex polyhedra. Most of the definitions presented here coincide with the corresponding definitions for rank $3$ abstract polytopes, although some or them have been slightly modified in order to be more accessible. Readers interested in general definitions and more details about the theory of abstract regular polytopes are referred to \cite[Section 2A]{ARP}.

An \emph{abstract polyhedron} is a partially ordered set $(\P,\leq)$ (we usually omit the order symbol) with a \emph{rank function} $rk$ with range $\{0,1,2\}$. In order to remind the geometrical origin of the theory, we call \emph{vertices}, \emph{edges} and \emph{faces} the elements of $\P$ of rank $0$, $1$, and $2$, respectively. We say that two elements $F$ and $G$ of $\P$ are \emph{incident} if $F \leq G$ or $G\leq F$. A \emph{flag} is a maximal chain of $\P$. Besides, $\P$ satisfies the following properties:
\begin{enumerate}[label=(P\arabic*)]
\item \label{flags} Every flag of $\P$ contains a vertex, an edge and a face.
\item \label{sfc} $\P$ is \emph{strongly flag-connected}, that is, given any two flags $\Phi$ and $\Psi$ of $\P$ there exists a sequence of flags $\Phi = \Phi_0 , \Phi_1 , \dots , \Phi_k = \Psi$ such that, for every $i \in  \{1, \dots , k\}$, $ \Phi_{i-1}$ and $\Phi_{i}$ are \emph{adjacent} (differ by precisely one element) and $\Phi \cap \Psi \subset \Phi_i$ for every $i \in \{0, 1, \dots , k\}$.
\item \label{diamond} $\P$ satisfies the \emph{diamond condition}, which means that every edge is incident to precisely two vertices and to two faces, and if a vertex $F_0$ is incident to a face $F_2$ then there are exactly two edges $F_1$ and $F'_1$ such that $F_0 \leq F_1 , F'_1 \leq F_2$.
\end{enumerate}

The $1$-skeleton of $\P$ is the (connected) graph $\sk(\P)$ determined by the vertices and edges of $\P$.

We sometimes identify a face $F_2$ with the set $\{G \in \P : G < F_2\}$. We can also think a face  $F_2$ as a graph whose nodes are precisely the vertices of $\P$ incident to $F_2$ and two of them are adjacent if they are incident to a common edge of $\P$. In this sense, we may think of faces of $\P$ as subgraphs of $\sk(\P)$.

If $F_0$ is a vertex of $\P$, the \emph{vertex figure} at $F_0$ is the set $\{G \in \P : F_0 < G \}$. We can also give the structure of a graph to the vertex figure of $F_0$ if we think the edges of $\P$ incident to $F_{0}$ as nodes and two nodes are adjacent if and only if the corresponding edges are incident to a common face of $\P$.

With the observations made above we may also define a polyhedron $\P$ as a connected graph $\sk(\P)$ and a family of subgraphs $\P_{2}$  that satisfy 

\begin{enumerate}[label=(P'\arabic*)]
\item \label{alledgesCovered} Every edge of $\sk(\P)$ belongs to at least one subgraph of $\P_{2}$.
\item Every graph of $\P_{2}$ is a connected $2$-valent graph.
\item Every vertex figure of $\P$ (defined in the analogous way as before) is a connected $2$-valent graph. 
\end{enumerate}

In this definition, the graphs of $\P_{2}$ are the \emph{faces} of $\P$. It is not hard to verify that both definitions are equivalent. Even though the first definition of polyhedra presented here is more popular in the literature we introduce the second one since it will be useful in Section\nobreakspace \ref {sec:RPT}.

Given a flag $\Phi$ of $\P$ and $0 \leq i \leq 2$ we define $\Phi^{i}$, the \emph{$i$-adjacent flag of $\Phi$}, as the (unique) flag that differs of $\Phi$ in precisely the element of rank $i$. If $F_0$ is a vertex, we define the \emph{degree} of $F_0$ as the number of edges incident to $F_0$. The \emph{codegree} of a face $F_2$ is the number of edges incident to $F_2$.

The \emph{group of automorphisms} of $\P$, denoted by $\G(\P)$, is the group of all order-preserving bijections $\gamma: \P \to \P$. Observe that $\G(\P)$ acts on $\P_i$, the set of elements of rank $i$, as well as on $\F(\P)$, the set of flags of $\P$. It is not hard to see that if $\gamma \in \G(\P)$ and $\Phi \in \F(\P)$, then $(\Phi^{i}) \gamma = (\Phi \gamma)^{i}$. It follows from the previous observation and from the strong flag-connectivity that the action of $\G(\P)$ on $\F(\P)$ is free. We say that a polyhedron  $\P$ is \emph{regular} if $\G(\P)$ acts transitively on $\F(\P)$. 

Observe that if $\P$ is a regular polyhedron and $\Phi$ is a flag of $\P $ there must exist \emph{distinguished automorphisms} $\r_0$, $\r_1$ and $\r_2$ such that $\Phi \r_i = \Phi^{i}$ for $i \in \{0,1,2\}$. The converse is also true, if $\P$ is a polyhedron and there exists a flag $\Phi$ and automorphisms $\r_0$, $\r_1$ and $\r_2$ such that $\Phi \r_i = \Phi^{i}$ for $i \in \{0,1,2\}$ then $\P$ is a regular polyhedron.  Notice that if $\P$ is regular, then every face has the same (possibly infinite) codegree and every vertex has the same (possibly infinite) degree. Let $p$ be the codegree of any face and $q$ the degree of any vertex, then we say that $\P$ has \emph{Schläfli type} $\{p,q\}$. In such situation, the distinguished automorphisms satisfy the relations 
\begin{equation} \label{eq:relations}
\r_{0}^{2}=\r_{1}^{2}=\r_{2}^{2}=(\r_0 \r_1)^{p} = (\r_1 \r_2)^q = (\r_0 \r_2)^2 = \epsilon.                                                                                                                                                                                                                                                                                                                                                                                                                                                                                                                                                                                                                                                                                                                                                                                                                                                                                                                                                                                                                                                                                                                                                                                                                                                                                                                                                                                                                                                                                                                                                                                                            \end{equation}

It is important to remark that Schläfli type may be also defined for non-regular polyhedra (see \cite{pellWeissChiralR3,schulteChiralR31}).

If $\P$ is a regular polyhedron with Schläfli type $\{p,q\}$, then $\G(\P)$ is a smooth quotient of the string Coxeter group $\Coxpq$ defined by \[\Coxpq:= \left\langle r_0, r_1, r_2 :r_{0}^{2}=r_{1}^{2}=r_{2}^{2}=(r_0 r_1)^{p} = (r_1 r_2)^q = (r_0 r_2)^2 = \epsilon \right\rangle.\]

In the case of the tetrahedron, the octahedron and the icosahedron, their automorphism groups are precisely the groups $\Tet$, $\Oct$ and $\Ico$ respectively (\cite[Theorem 3B3]{ARP}). 

Given the Coxeter group $\Coxpq$, we denote by $\coxpq$ the subgroup of all elements of $\Coxpq$ that can be written as a product of an even number of generators. If $\Coxpq$ is represented by a group generated by reflections (as the automorphism group of a regular convex polyhedron), $\coxpq$ is the \emph{rotational subgroup}. With this notation it is easy to prove that $\tet \leq \Tet \leq \Oct$.

\subsection{Realizations of Polyhedra in Euclidean 3-space}

In \cite[Section 5A]{ARP} the authors define a \emph{(Euclidean) realization} of a regular polyhedron $\P$ as a function $\beta: \P_0 \to E$ where $\P_i$ is the set of elements of $\P$ whose rank is $i$ and $E$ is some Euclidean space. We can recover the structure of the polyhedron defining $\beta = \beta_0$, $V_0:= \P_0 \beta$, and recursively for $i \in \{1,2\}$, defining $\beta_i:\P_i \to V_i$, for some $V_i \subset 2^{V_{i-1}}$, as \[F \beta_{i} = \{G \beta_{i-1}: G \in \P_{i-1} \ \text{y} \ G < F \}.\] Informally speaking, every edge $F_1$ is the pair of vertices incident to $F_1$ and every face $F_2$ is the family of edges incident to $F_2$. However, in general an edge (face) is not determined by the vertices (edges) incident to the edge (face), so for our purposes it is convenient to think an edge as a segment joining the two vertices incident to the edge and a face as the family of segments determined by edges incident to the given face. Somehow, this codifies the idea of realization of a polyhedron $\P$ in terms of the graph $\sk(\P)$ in the same sense that we may define an abstract polyhedron in terms of this graph. Observe that for Euclidean spaces this convention is equivalent to the definition presented before. 

In this work, when talking about realizations of polyhedra we usually identify the vertices of a polyhedron with the corresponding points of the realization, the edges with the line segments mentioned before and the faces with the collection the corresponding line segments. When there is no confusion we also think a polyhedron both, abstractly and geometrically realized and we just refer to it as a polyhedron.

A realization $\beta$ is faithful if each $\beta_i$ is injective. The realization is \emph{discrete} if $V_0$ is a discrete subset of $E$. If $\beta$ is a realization of $\P$, define the \emph{symmetry group} of $\P$, denoted by $G(\P)$, as the group of all isometries of $E$ preserving the structure of $\P$. If we assume that $E = \aff(\P)$ then $G(\P)$ is an Euclidean representation of a subgroup $\Delta \leq \G(\P)$. We say that $\beta$ is \emph{symmetric} whenever $\Delta = \Gamma(\P)$.

Notice that if $\beta$ is symmetric and faithful then $\G(\P) \cong   G(\P)$. If $\P$ is a symmetrically realized regular polyhedron with distinguished generators $\r_0$, $\r_1$ and $\r_2$, we denote by $R_0$, $R_1$ and $R_2$ the corresponding symmetries. Unless otherwise specified, in this section we will only talk about symmetric realizations and we will refer to them only as realizations.

A realization is said to be \emph{blended} if there are proper orthogonal complementary subspaces $L$ and $M$ of $E$ such that $G(\P)$ permutes the orthogonal translates of $L$ (and hence, of $M$). A blended realization $\beta$ induces realizations $P_1$ and $P_2$ of $\P$ on $L$ and $M$ respectively, in this case we write $\P = P_1 \# P_2$. If a realization is not blended we call the polyhedron \emph{pure}.  

In \cite{ordinary} P. McMullen and E. Schulte list $48$ faithful and discrete realizations of regular polyhedra in Euclidean $3$-space and prove that such list is complete. The list coincides with the $48$ regular polyhedra found by B. Grünbaum and A. Dress in \cite{grunbaumOldAndNew}, \cite{dress1} and \cite{dress2}. In Section\nobreakspace \ref {sec:RPT} we extend the definition of realization in order to talk about realizations in the $3$-Torus $\TT$ (defined below) and give a classification of regular polyhedra in $\TT$. 

Before talking about $\TT$ we give a brief review of realizations of regular polyhedra in $\E$. Readers interested in a more detailed discussion about this topic are referred to \cite{grunbaumOldAndNew} and \cite{ordinary}.

\subsubsection{Finite Polyhedra}

Among the $48$ regular polyhedra in $\E$, $18$ of them are finite. The list contains the $5$ \emph{Platonic solids}; the $4$ \emph{Kepler-Poinsot polyhedra} which are related to the icosahedron and the dodecahedron by the \emph{$2$-facetting operation} and duality. These polyhedra and the facetting operation are described in detail in \cite[Chapter VI]{coxeterRegularPolytopes}. The remaining $9$ are the \emph{Petrie-duals (Petrials)} of the nine listed above. The Petrial of $\P$ has the same $1$-skeleton of $\P$ but take as faces the \emph{Petrie polygons}, which are the edge-paths of $\sk(\P)$ defined by the property that any two but not three consecutive edges belong to the same face. For an algebraic description of the $2$-facetting operation and the Petrie-operation see \cite[Section 7B]{ARP}.

In \cite{ordinary} and \cite[Section 7E]{ARP} the authors classify the 18 finite regular polyhedra according to their symmetry group and list them in the diagrams \eqref{eq:tetrasym}, \eqref{eq:octasym} and \eqref{eq:icosym}. In such diagrams the arrows represent that the polyhedra are related by the Petrie-operation ($\pi$), the $2$-facetting operation ($\varphi_2$) or the well known \emph{duality operation} ($\delta$).

Tetrahedral symmetry:
\begin{equation} \label{eq:tetrasym}
	\begin{tikzcd}
		\{3,3\} \arrow[leftrightarrow]{r}{\pi}  & \{4,3\}_3
	\end{tikzcd}
\end{equation}

Octahedral symmetry:
\begin{equation}\label{eq:octasym}
\begin{tikzcd}
  	\{6,4\}_{3} \arrow[leftrightarrow]{r}{\pi} & \{3,4\}\arrow[leftrightarrow]{r}{\delta} & \{4,3\}\arrow[leftrightarrow]{r}{\pi} & \{6,3\}_4
  \end{tikzcd}
\end{equation}

Icosahedral symmetry:
\begin{equation} \label{eq:icosym}
\begin{tikzcd}
	\{10,5\} \arrow[leftrightarrow]{r}{\pi} \arrow[leftrightarrow]{d}{\varphi_2} & \{3,5\}\arrow[leftrightarrow]{r}{\delta} \arrow[leftrightarrow]{d}{\varphi_2} & \{5,3\}\arrow[leftrightarrow]{r}{\pi} & \{10,3\} \\
	\{6,\frac{5}{2}\}\arrow[leftrightarrow]{r}{\pi} & \{5,\frac{5}{2}\}\arrow[leftrightarrow]{r}{\delta} & \{\frac{5}{2},5\}\arrow[leftrightarrow]{r}{\pi} \arrow[leftrightarrow]{d}{\varphi_2} & \{6,5\} \arrow[leftrightarrow]{d}{\varphi_2}\\
	\{\frac{10}{3},3\} \arrow[leftrightarrow]{r}{\pi} & \{\frac{5}{2},3\}\arrow[leftrightarrow]{r}{\delta} & \{3,\frac{5}{2}\}\arrow[leftrightarrow]{r}{\pi} & \{\frac{10}{3},\frac{5}{2}\}
\end{tikzcd}
\end{equation}

For purposes of this work, we give in \MakeUppercase Table\nobreakspace \ref {tab:finite} explicit coordinates and generators of $G(\P)$ for $\P$ being the tetrahedron $\{3,3\}$, the cube $\{3,4\}$ and the octahedron $\{4,3\}$. 

\begin{table}
\begin{center}
\resizebox{\textwidth}{!} {
\begin{tabular}{|c| c| c |c |c| c|}
	\hline 
	Polyhedron $\P$ & Vertex set & Base vertex & \multicolumn{3}{c|}{Generators of $G(\P)$} \\ 
	&  &  & \multicolumn{1}{c}{$R_0:(x,y,z) \mapsto$} & \multicolumn{1}{c}{$R_1:(x,y,z) \mapsto$} & \multicolumn{1}{c|}{$R_2:(x,y,z) \mapsto$} \\ 
	\hline 
	$\{3,3\}$ & $\!\begin{aligned}[t] &\left\{(1,1,1), (-1,-1,1), \right. \\ &\left. \phantom{ \{ } (-1,1,-1), (1,-1,-1) \right\}  \end{aligned}$ & $(1,1,1)$ & $(-y,-x,z)$ & $(z,y,x)$ & $(y,x,z)$ \\ 
	$\{3,4\}$ & $\{\pm(1,0,0), \pm (0,1,0), \pm (0,0,1)\}$ & $(1,0,0)$ & $(y,x,z)$ & $(x,z,y)$ & $(x, y, -z)$ \\ 
	$\{4,3\}$  & $\left\{(\pm 1, \pm 1, \pm 1) \right\}$ & $(1,1,1)$ &  $(x, y, -z)$ & $(x,z,y)$ & $(y,x,z)$  \\ \hline 
\end{tabular}
}
\caption{Vertex-set and group of symmetries of $\{3,3\}$, $\{3,4\}$ and $\{4,3\}$} \label{tab:finite}
\end{center}
\end{table}

\subsubsection{Planar Apeirohedra}

There are $6$ planar infinite polyhedra (apeirohedra). They are the planar tessellations with squares (which is selfdual), equilateral triangles and regular hexagons (which are dual of each other). The remaining $6$ are just the Petrials of these tessellations. \MakeUppercase Table\nobreakspace \ref {tab:planar} shows the vertex set and distinguished generators for $G(\P)$ for each planar polyhedra. In this table $X$ denotes the set $\left\lbrace a(1,0) + b(\frac{1}{2} , \frac{\sqrt{3}}{2}) : a,b \in \mathbb{Z}\right\rbrace$.

\begin{table}
\begin{center}
\resizebox{\textwidth}{!} {
\begin{tabular}{|c| c| c |c |c| c|}
\hline 
Polyhedron $\P$ & Vertex set & Base vertex & \multicolumn{3}{c|}{Generators of $G(\P)$} \\ 
&  &  & \multicolumn{1}{c}{$R_0:(x,y,z) \mapsto$} & \multicolumn{1}{c}{$R_1:(x,y,z) \mapsto$} & \multicolumn{1}{c|}{$R_2:(x,y,z) \mapsto$} \\ 
\hline
$\{4,4\}$ & $\mathbb{Z}^2$ & $o$ & $(1-x,y)$ & $(y,x)$ & $(x,-y)$ \\ 
$\{3,6\}$  & $X$ & $o$ &  $(1-x,y)$ & $\left(\frac{1}{2} x+ \frac{\sqrt{3}}{2}y, \frac{\sqrt{3}}{2}x -  \frac{1}{2}y \right)$ & $(x,-y)$  \\ 
$\{6,3\}$ & $\left((\frac{1}{2}, \frac{\sqrt{3}}{6})+X \right) \cup \left((1,\frac{\sqrt{3}}{3})+X \right) $ & $(\frac{1}{2}, \frac{\sqrt{3}}{6})$ & $(x,-y)$ & $\left(\frac{1}{2} x+ \frac{\sqrt{3}}{2}y, \frac{\sqrt{3}}{2}x -  \frac{1}{2}y \right)$ & $(1-x,y)$ \\
\hline 
\end{tabular}
}
\caption{Vertex-set and group of symmetries of Regular Plane Apeirohedra} \label{tab:planar}
\end{center}
\end{table}

\subsubsection{Blended Apeirohedra}

There are $3$ blended polyhedra with finite faces. Each of these polyhedra are obtained by the blend of a planar tessellation $\Q$ and a line segment $\{\ \}$. 

If $\Q=\{4,4\}$ (resp. $\{6,3\}$) then $\sk(\Q)$ is a bipartite graph and $\Q \# \{\ \}$ is constructed by deforming planar squares (resp. hexagons) of $\Q$ to skew squares (resp. hexagons) by lifting orthogonally alternate vertices of $\Q$ to equal height $\alpha$ above the plane of $\Q$. If $\Q = \{3,6\}$, then the vertices of   $\Q \# \{\ \}$ are two parallel copies of the vertices of $\Q$ in planes at distance $\alpha$; the faces are skew hexagons that go twice around each triangular prism determined by two parallel faces of the copies of $\Q$. In any of these polyhedra, the parameter $\alpha$ varies continuously in the interval $(0,+\infty)$. 

Three more blended polyhedra of type $\Q \# \{\infty\}$, with $\Q$ a planar tessellation, can be constructed by lifting a tower of prisms of height $\alpha$ over each face of $\Q$. The faces of $\Q \# \{\infty\}$ are helices that that go around each tower of prisms. If $\Q$ is $\{4,4\}$ or $\{3,6\}$ every tower of prisms has one helix surrounding it, and helices that project to adjacent faces of $\Q$ go in opposite sense. If $\Q = \{6,3\}$, every tower of prisms has six helices (three left and three right). The parameter $\alpha$ varies continuously in the interval $(0,+\infty)$.

The remaining $6$ blended polyhedra are just the Petrie-duals of those described above. 

Table\nobreakspace \ref {tab:blended} gives explicitly distinguished generators for the group of symmetries of the blended polyhedra described above.

\begin{table}
\begin{center}
\resizebox{\textwidth}{!} {
\begin{tabular}{|c| c |c |c| c|}
	\hline 
Polyhedron $\P$ &Base vertex & \multicolumn{3}{c|}{Generators of $G(\P)$} \\ 
&  & \multicolumn{1}{c}{$R_0:(x,y,z) \mapsto$} & \multicolumn{1}{c}{$R_1:(x,y,z) \mapsto$} & \multicolumn{1}{c|}{$R_2:(x,y,z) \mapsto$} \\ 
\hline
$\{4,4\} \# \{\ \}$ & $o$ & $(1-x,y,\alpha - z)$ & $(y,x,z)$ & $(x,-y,z)$ \\ 
$\{6,3\}\# \{\ \}$  & $\left(\frac{1}{2}, \frac{\sqrt{3}}{6}, 0 \right)$ & $(x,-y, \alpha - z)$ & $\left(\frac{1}{2} x+ \frac{\sqrt{3}}{2}y, \frac{\sqrt{3}}{2}x -  \frac{1}{2}y , z \right)$ & $(1-x,y,z)$\\
$\{3,6\}\# \{\ \}$ & $o$ & $(1-x,y, \alpha-z)$ & $\left(\frac{1}{2} x+ \frac{\sqrt{3}}{2}y, \frac{\sqrt{3}}{2}x -  \frac{1}{2}y , z \right)$ & $(x,-y,z)$ \\
$\{4,4\} \# \{\infty \}$ & $o$ & $(1-x,y,\alpha-z)$&$(y,x,-z)$ &$(x-y,z)$\\
$\{6,3\} \# \{\infty \}$  & $\left(\frac{1}{2}, \frac{\sqrt{3}}{6},0 \right)$ & $(x,-y, \alpha - z)$ & $\left(\frac{1}{2} x+ \frac{\sqrt{3}}{2}y, \frac{\sqrt{3}}{2}x -  \frac{1}{2}y , -z \right)$ & $(1-x,y,z)$\\
$\{3,6\} \# \{\infty \}$ & $o$ & $(1-x,y, \alpha-z)$ & $\left(\frac{1}{2} x+ \frac{\sqrt{3}}{2}y, \frac{\sqrt{3}}{2}x -  \frac{1}{2}y , -z \right)$ & $(x,-y,z)$ \\
\hline 
\end{tabular}
}
\caption{Generators and group of symmetries of regular blended apeirohedra} \label{tab:blended}
\end{center}
\end{table}

\subsubsection{Pure Apeirohedra}

Infinite pure polyhedra (or pure apeirohedra) are those regular polyhedra that do not admit a blended realization. There are $12$ of them which we describe below.

\emph{Petrie-Coxeter polyhedra} $\{4,6|4\}$, $\{6,4|4\}$ and $\{6,6|3\}$ are three polyhedra found by  by H.S.M. Coxeter and J.F. Petrie and described in detail in \cite{coxeterSkew}. The polyhedron $\{4,6|4\}$ consists of the vertices and edges of the cubical tessellation $\{4,3,4\}$; its faces are half of the squares of the tessellation taken in an alternating way in each plane parallel to the coordinated planes. The polyhedron $\{6,4|6\}$ consists of the vertices, edges and hexagons of each truncated octahedron of the uniform tessellation $(4.6.4)^{4}$ with truncated octahedra. The polyhedron $\{6,6|3\}$ consists of the vertices, edges and hexagons of each truncated tetrahedra of the uniform tessellation of the space $(3.3.3)^{2}.(3.6.6)^{6}$ with regular tetrahedra and truncated tetrahedra. The tessellations $(4.6.4)^{4}$ and $(3.3.3)^{2}.(3.6.6)^{6}$ are described in detail in \cite{grunbaumUniformTilings}. 

Applying the Petrie operation to $\{4,6|4\}$, $\{6,4|4\}$ and $\{6,6|3\}$ we obtain the polyhedra $\{\infty,6\}_{4,4}$, $\{\infty,4\}_{6,4}$ and $\{\infty,6\}_{6,3}$, respectively.

The other $6$ pure apeirohedra can be paired according to their $1$-skeleton.

The vertices of the polyhedron $\{6,6\}_4$ are the vertices of the uniform tessellation $(3.3.3)^{8}.(3.3.3.3)^{6}-A$ with regular octahedra and regular tetrahedra (see \cite{grunbaumUniformTilings}). The facets are given by taking one Petrie polygon of each octahedron. The polyhedron $\{4,6\}_6$ is the Petrial of $\{6,6\}_4$ and so, they share the $1$-skeleton. The faces of $\{4,6\}_6$ are Petrie polygons of the tetrahedra in the uniform tessellation $(3.3.3)^{8}.(3.3.3.3)^{6}-A$, one of each tetrahedron.

The polyhedron $\{\infty,3\}^{(b)}$ is described with full detail in \cite{ordinary}. Its faces are helices over squares whose axes are parallel to the coordinated axes. The polyhedron $\{\infty,3\}^{(a)}$ is the Petrial of $\{\infty,3\}^{(b)}$; its faces are triangular helices. Both polyhedra are also described in \cite{pellWeissChiralR3}.

The polyhedron $\{6,4\}_6$ is the dual of $\{4,6\}_6$ and its $1$-skeleton is the translate by $(\frac{1}{2},\frac{1}{2},\frac{1}{2})$ of subgraph of the cubic tessellation induced by the points $\bZ^3 \sm X$ where $X$ is the set of points of integer coordinates of the form $(2m, 2n, 2k+1)$ y $(2m+1, 2n+1, 2k)$. The face-set of $\{6,4\}_6$ consists of one Petrie polygon of half of the cubes whose centre is a point of $\mathbb{Z}^3$. This polyhedron is \emph{self-Petrie}, that is, isomorphic to its Petrial.

The polyhedron $\{\infty, 4\}_{\cdot , \ast 3}$ has the same $1$-skeleton than $\{6,4\}_6$. Its facets are helices over triangles whose axes are parallel to the vectors $(\pm 1, \pm 1, \pm 1)$. This polyhedron is also self-Petrie. 

In Table\nobreakspace \ref {tab:pure} we give explicit generators for the group of symmetries of every pure apeirohedra described before.

\begin{table}
\begin{center}
\resizebox{\textwidth}{!} {
\begin{tabular}{|c | c | c | c| c|}
\hline 
Polyhedron $\P$ & Base Vertex & \multicolumn{3}{c|}{Generators of $G(\P)$} \\ 
 &  & \multicolumn{1}{c}{$R_0:(x,y,z) \mapsto$} & \multicolumn{1}{c}{$R_1:(x,y,z) \mapsto$} & \multicolumn{1}{c|}{$R_2:(x,y,z) \mapsto$} \\ 
\hline
$\{4,6|4\}$ & $o$ & $(1-x,y,z)$ & $(y,x,-z)$ & $(x,z,y)$ \\ 
$\{6,4|4\}$ & $\left(\frac{1}{2}, \frac{1}{2}, 0\right)$& $(x,z,y)$ & $(y,x,-z)$ & $(1-x,y,z)$\\ 
$\{6,6|3\}$ & $o$ & $(x,1-z,1-y)$ & $(y,x,-z)$ & $(x,z,y)$ \\ 
$\{\infty,6\}_{4,4}$ & $o$ & $(1-x,z,y)$ & $(y,x,-z)$ & $(x,z,y)$ \\ 
$\{\infty, 4\}_{6,4}$ & $\left(\frac{1}{2}, \frac{1}{2}, 0\right)$ & $(1-x,z,y)$ & $(y,x,-z)$ & $(1-x,y,z)$ \\ 
$\{\infty,6\}_{6,3}$ & $o$ & $(x,1-z,1-y)$ & $(y,x,-z)$ & $(x,z,y)$ \\ 
$\{6,6\}_{4}$ & $o$ & $(1-y,1-x,-z)$ & $(x,z,y)$ & $(y,x,-z)$ \\ 
$\{4,6\}_{6}$ & $o$ & $(1-x,1-y,z)$ & $(x,z,y)$ & $(y,x,-z)$ \\ 
$\{\infty,3\}^{(b)}$ & $o$ & $(1-y,1-x,-z)$ & $(z,-y,x)$ & $(y,x,-z)$ \\ 
$\{\infty,3\}^{(a)}$ & $o$ & $(1-x,1-y,z)$ & $(z,-y,x)$ & $(y,x,-z)$ \\ 
$\{6,4\}_{6}$ & $\left(\frac{1}{2}, \frac{1}{2}, \frac{1}{2}\right)$ & $(y,x,-z)$ & $(x,z,y)$ & $(1-x,1-y,z)$ \\ 
$\{\infty, 4\}_{\cdot, \ast 3}$ & $\left(\frac{1}{2}, \frac{1}{2}, \frac{1}{2}\right)$ & $(y,x,-z)$ & $(1-x,z,y)$ & $(1-x,1-y,z)$ \\ \hline 
\end{tabular}
}
\caption{Generators and group of symmetries of regular pure apeirohedra} \label{tab:pure}
\end{center}
\end{table}

\section{Lattice groups and the $3$-torus $\TT$.} \label{sec:lattices}

Given $1 \leq d \leq 3$ and $V \leq \E$ a $d$-dimensional vectorial space, a \emph{lattice group} (of rank $d$) $\bfL$ on $V$ is a subgroup of the translation group $\T(\E)$  of $\E$ generated by $d$ linearly independent translations $t_1, \dots t_d$ such that each translation vector lies on $V$. If $\bfL$ is a lattice group and $o$ is the origin of $\E$, the \emph{lattice} $\L$ associated to $\bfL$ is the subset $o\bfL$. If $v_i$ is the translation vector of $t_i$, the set $\{v_1, \dots , v_d\}$ is called a \emph{basis} of $\L$. Notice that $\{v_1, \dots , v_d\}$ is also a basis of $V$ and $\L$ consists of the integer linear combinations of $\{v_1, \dots , v_d\}$. Finally, observe that $\L$ is a discrete subset of $\E$. 

An open subset $D$ of $\E$ is an \emph{fundamental region} for a lattice group $\bfL$ if $\E = \bigcup_{t \in \bfL} (\clos D) t$ and $Dt \cap  Dt' = \emptyset$ if $t \neq t'$, where $\clos D$ denotes the closure of $D$ . Given $\bfL$ a lattice group on $V$ with base $\{v_1, \dots ,v_d\}$, the open parallelepiped $P(v_1, \dots , v_d)$ determined by de vectors $v_1, \dots v_d$ is a fundamental region for $\bfL$. Another important fundamental region is the so called \emph{ (open) Dirichlet domain} (centred at the origin $o$) \[D(\bfL):= \{x \in V: d_{\E}(o,x) < d_{\E}(o, x t) : t \in \bfL \sm \{id\} \}.\]

Following \cite[Section 6D]{ARP}, if $\{e_1, \dots , e_n\}$ denotes the standard basis of $\bb{E}^n$, we denote by $\sql$ the \emph{$2$-dimensional square lattice} $\bb{Z}^2$ generated by $\{e_1,e_2\}$, and by $\scl$ the \emph{square-centred lattice}, generated by $\{e_1+e_2, e_1-e_2\},$ that consist of all integral vectors whose coordinates have the same parity; observe that $\sql$ and $\scl$ are \emph{similar}, that is, differ by a similarity, namely, the composition of a dilatation by factor $\sqrt{2}$ and a rotation of $\frac{\pi}{4}$ . Similarly, we denote by $\cl$ the \emph{$3$-dimensional cubic lattice} $\bb{Z}^3$ generated by $\{e_1, e_2, e_3\}$; $\bcl$ is the \emph{body-centred cubic lattice} with basis $\{2e_1, 2e_2, e_1+e_2+e_3\}$ that consists of all integral vectors whose coordinates have the same parity; and finally $\fcl$ is the \emph{face-centred cubic lattice}, generated by $\{e_1+e_2 , e_1 - e_2, e_3-e_2\}$ and consists of all integral vectors whose coordinate sum is even. We denote by $\Sql$, $\Scl$, $\Cl$, $\Bcl$ and $\Fcl$ the corresponding translation groups. 

We commonly identify the lattice $\sql$ with the sublattice of $\cl$ generated by $\{e_1,e_2\}$, and write $\L_{(1,0)} \subset \cl$ and $ \scl \subset \fcl$. Moreover, if $\bfL$ is a lattice group over $V$ and $w$ is a vector not in $V$, the lattice $\L \+ w$ is the lattice associated to the group $\left\langle \bfL,  t_w \right\rangle$ with $t_w$ the translation by $w$. Then we can think $\cl$ as the lattice $\L_{(1,0)} \+ e_3$. Finally, if $T$ is a linear transformation that preserves $V$, we denote $\L T$ the lattice $\{vT : v \in \L \}$ and by $\bfL^T$ the corresponding lattice group. In the particular case that $T$ is the dilatation by a factor $a \in \mathbb{R}$, we denote by $a\L$ the lattice $\L T$ and the corresponding lattice group by $a\bfL$.

Unless specified otherwise, from now on $\bfL$ denotes a rank-$3$-lattice group on $\E$ and $\L$ its corresponding lattice. The \emph{$3$-Torus} associated to $\bfL$, denoted by $\TTL$ (or simply $\TT$ if there is not confusion regarding $\bfL$), is the  quotient space $\E / \bfL$ whose points are orbits of $\E$ under $\bfL$, provided with the metric $d_{\bfL}$, which is induced by the Euclidean metric $d$ in $\E$ and defined by \[d_{\bfL}(x\bfL, y\bfL) = \inf \{d(xt , yt') :  t,t' \in \bfL\}.\] If $D$ is a fundamental region (for instance, the Dirichlet domain) for $\bfL$, $\TTL$ is isometric to $D \cup (\bd D/\sim)$ where $\bd D$ is the boundary of $D$ and $x \sim y$ if and only if $x=yt$ for some $t \in \bfL$.

In order to study the group of isometries of $\TTL$, a first approach is to study those isometries of $\E$ that induce an isometry in $\TTL$. In other words, we want to classify the isometries $S$ of $\E$ that make commutative the diagram in \eqref{eq:diagtorus}, where $\pi_{\bfL}$ denotes the quotient map. It is easy to prove that those isometries are precisely those that normalize $\bfL$ in $\IE$. Moreover, the following theorem holds (see \cite[p.336]{ratcliffeHypManifolds}).

\begin{thm}\label{thm:isomTorus}
	The group of isometries of $\TTL$ is isomorphic to $N_{\I(\E)}(\bfL) / \bfL$ where $N_{\I(\E)}(\bfL)$ denotes the normalizer of $\bfL$ in $\I(\E)$. Furthermore, the isomorphism is induced by the diagram \eqref{eq:diagtorus}.
\end{thm}

\begin{equation} \label{eq:diagtorus}
\begin{tikzcd}
\E \arrow{r}{S} \arrow{d}{\pi_{\bfL}} & \E \arrow{d}{\pi_{\bfL}} \\
\TTL \arrow[dashed]{r}{\overline{S}} & \TTL
\end{tikzcd}
\end{equation}

Given the results mentioned so far, we are interested in studying those isometries of $\E$ that normalize $\bfL$. Observe that every translation induces an isometry in $\TTL$, since the translation group is abelian. Recall that $\I(\E) = \T(\E) \rtimes \O(3)$, where $\O(3)$ is the group of orthogonal isometries of $\E$, and hence every isometry $S$ of $\E$ can be (uniquely) written as a product $t S'$ with $t \in \T(\E)$ and $S' \in \O(3)$. Since $\T(\E) \leq N_{\I(\E)}(\bfL)$ for every lattice group $\bfL$, an isometry $S=tS'$, with $t \in \T(\E)$ and $S' \in \O(3)$, belongs to $N_{\I(\E)}(\bfL)$ if and only if $S' \in N_{\I(\E)}(\bfL)$. Hence we can restrict our work to linear isometries.  Notice that if $t_v$ is a translation by the vector $v$ and $S$ is a linear isometry of $\E$, then $S^{-1}t_v S = t_{vS}$ where $t_{vS}$ is the translation by $vS$. Therefore a linear isometry $S$ normalizes $\bfL$ if and only if $S$ preserves $\L$.

Now it is immediate that every rotatory reflection of period $2$ in $\E$ induces an isometry of $\TTL$ for every $\bfL$ since its orthogonal component is the central inversion $-id$, which always preserves the lattice $\L$. 

We find it convenient to introduce some notation now. If $G \leq \IE$ is a group of isometries of $\E$, let $G_{o}$ denote the \emph{special group} associated to $G$, which is the group of all orthogonal components of isometries of $G$. By $\bar{G}_{o}$ we denote de \emph{extended special group} associated to $G$, which is the group $\left\langle G_{o}, -id \right\rangle$.  

We can summarize these observations in the following results. 

\begin{lem} \label{lem:equivNorm}
Let $S=tS'$ an isometry of $\E$ with $t \in \T(\E)$ and $S' \in \O(3)$. Let $\bfL$ be a lattice group. The following statements are equivalent:
\begin{enumerate}[label=(\roman*)]
\item $S$ normalizes $\bfL$.
\item $S$ induces an isometry $\bar{S}$ of $\TTL$ that makes commutative the diagram of \eqref{eq:diagtorus}.
\item $S'$ normalizes $\bfL$.
\item $S'$ induces an isometry of $\TTL$ that makes commutative the diagram of \eqref{eq:diagtorus}.
\item $S'$ preserves $\L$.
\end{enumerate}
\end{lem}

Since we are interested in the study of groups of isometries of regular polyhedra in the $3$-torus, we consider more appropriate to give a result equivalent to the lemma above but in terms of groups of isometries.

\begin{lem} \label{lem:equivNormGroup}
Let $G \leq \IE$ be a group of isometries of $\E$. Let $G_{o}$ and $\bar{G}_{o}$ denote the special group and the extended special group associated to $G$, respectively. Let $\bfL$ be a lattice group and $N_{\I(\E)}(\bfL)$ its normalizer in $\IE$. The following statements are equivalent:
\begin{enumerate}[label=(\roman*)]
\item $G \leq N_{\I(\E)}(\bfL)$.
\item $G_{o} \leq N_{\I(\E)}(\bfL)$.
\item $\bar{G}_{o} \leq N_{\I(\E)}(\bfL)$.
\item $G$ admits a representation as isometries of $\TTL$.
\end{enumerate}
\end{lem}

Since the group of symmetries of a regular polyhedron is generated by involutions, we are interested in studying those lattices invariant under involutory orthogonal isometries. However, since $-id$ preserves every lattice we just have to study those lattices preserved by plane reflections and half-turns. The following lemma states that it is enough to understand lattices preserved by reflections. 

\begin{lem}\label{lem:RefequivRot}
Let $R$ be the reflection with respect to a plane $\Pi$ that contains the origin $o$ of $\E$ and let $\L$ be a lattice. If $S$ is the half-turn with respect to $\Pi^{\perp}$, then $\L$ is preserved by $R$ if and only if $\L$ is preserved by $S$.
\begin{proof}
 It follows from the facts that $S=-R:=-id \circ R$ and $-id \in \norm$ for every lattice group $\bfL$.
\end{proof}

\end{lem} 

As a consequence of the lemma above, it is enough to study lattices invariant under the action of reflections. In \cite{equivelar4toroids}, Hubard, Orbanić, Pellicer and Weiss discuss widely the structure of such lattices, in particular they prove the \MakeUppercase Lemma\nobreakspace \ref {lem:ReflectionLattices}. 

Notice that if $\L$ is a rank-$3$ lattice on $\E$ that is invariant under the reflection $R$ in $\Pi$, then $\L_0:= \L\cap \Pi$ is a rank two lattice. 

Some lattices $\L$ that are trivially preserved by the reflection in $\Pi$ are those of the form $\L_0 \+ w$ where $\L_0 = \L \cap \Pi$ and $w$ is a vector orthogonal to $\Pi$. We call those lattices \emph{vertical translation lattice with respect to $\Pi$} (or simply \emph{vertical translation lattices} if there is not confusion regarding $\Pi$). We use similar terminology for rank-$2$ lattices preserved by reflections in lines.
 
The following lemma gives some useful structure of those lattices preserved by a reflection. 

\begin{lem} \label{lem:ReflectionLattices}
Let $\Pi$ be a plane and $\L$ a lattice invariant under the reflection in $\Pi$. Let $\L_0 := \L \cap \Pi$ and $w$ be a point fo $\L \sm \L_0$ such that $d(w,\Pi) \leq d(y,\Pi)$ for every point $y \in \L \sm \L_0$. Then \[\L = \L_0 \+ w = \bigcup_{k \in \mathbb{Z}} (\L_0 + kw).\]  Furthermore, if $\{v_1,v_2\}$ is any base for $\L_0$, $w$ can be chosen such that its projection to $\Pi$ belongs to $\{o, \frac{v_1}{2}, \frac{v_2}{2}, \frac{v_1 + v_2}{2}$, and $\L$ is a vertical translation lattice with respect to $\Pi$ if and only if $w$ projects to $o$.
\end{lem}

With the notation of the previous lemma, if $\L$ is not a vertical translation lattice, the points in the \emph{layer} $\L_0 +kw$ project orthogonally to $\Pi$ in points of $\L$ if and only if $k$ is even. If $k$ is odd, the points in the layer $\L_0+kw$ project to midpoints of points of $\L_0$. Furthermore, $\L_0 + kw$ and $\L_0 + (k+2)w$ differ by a translation orthogonal to $\Pi$ for any $k \in \mathbb{Z}$. 

It is clear that a similar result hold in dimension $2$. By using this analogous we can prove the following classification result regarding rank-$2$ lattices preserved by a line-reflection. The result uses strongly the fact that every lattice of rank $1$ is the dilatation of the integer lattice.  

\begin{lem}\label{lem:bidimlattices}
Let $\bfL$ be a rank-$2$ lattice group on $\mathbb{E}^2$ and $\Pi$ the line $x=0$. Let $\L$ be the lattice associated to $\bfL$. If $\L$ is preserved by the reflection in $\Pi$ then there exists a diagonal matrix $D$ such that $\L D = \L_{(1,0)}$ or $\L D = \L_{(1,1)}$.
\begin{proof}
Let $\L_0 = \L \cap \Pi$. Since $\L_0$ is a rank-$1$ lattice there exists $d_1 \in \bR$ such that $d_{1} e_{1}$ is a basis of $\L_0$. Let $w \in \L \sm \L_0$, a point of minimal distance $d_2$ to $\Pi$. If $\L$ is a vertical translation lattice with respect to $\Pi$, then define $D= \diag(d_1,d_2)$ and $\L = \sql D$. Otherwise, $w$ can be chosen such tat $w$ projects to $\left( \frac{d_1}{2}, 0 \right)$. Then define $D=\diag\left(\frac{d_1}{2}, \frac{d_2}{2} \right)$ and we have $\L= \scl D$. 
\end{proof}
\end{lem}

To conclude this section we find it convenient to introduce notation for some $2$-dimensional lattices. The \emph{triangular lattice} $\tl$ is the lattice whose basis is the set $ \TB = \{e_1, \frac{1}{2} e_1 + \frac{\sqrt{3}}{2} e_2\}$. It consists of the vertices of the triangular tessellation $\{3,6\}$. The \emph{triangle-centred lattice} $\tcl$ is the sublattice of $\tl$ generated by the vectors $(1,1)$ and $(2,-1)$ with respect to the basis $\TB$, i.e., the vectors $\frac{3}{2} e_1 + \frac{\sqrt{3}}{2} e_2$ and $\frac{3}{2} e_1 - \frac{\sqrt{3}}{2} e_2$. As always, we denote the corresponding lattice groups as $\mathbf{\tl}$ and $\mathbf{\tcl}$.

 \section{Regular polyhedra in the $3$-torus} \label{sec:RPT}
 
In this secction we define the concept of \emph{toroidal realization} of an abstract polyhedron as a generalization of Euclidean realization treated in Section\nobreakspace \ref {sec:basics}. It is important to mention that this idea has been explored for other spaces, for instance, in \cite{roliProjective1} and \cite{roliProjective2} Javier Bracho et. al. study realizations of regular polyhedra in the projective space. Many of the definitions presented in this section are inspired by Bracho's work and by the theory presented in \cite[Section 5A]{ARP}.
 
A \emph{graph in the $3$-torus} $\cX$ (or a \emph{graph in $\TT$}, for short) is a set $V(\cX)$ of points of $\TT$ called \emph{vertices} together with a set of geodesic-arcs $E(\cX)$ called \emph{edges} whose endpoints are two (distinct) points of $V(\cX)$ and for any edge $e \in E(\cX)$, $\inte e \cap V(\cX) = \emptyset$. Observe that any graph in $\TT$ has an underlying combinatorial graph where the combinatorial incidence is the same as the geometric incidence. 

Given two combinatorial graphs $X$ and $Y$, a \emph{graph homomorphism} $f: X \to Y$ is an adjacency-preserving function from de vertices of $X$ to the vertices of $Y$. A graph-homomorphism $f: X \to Y$ is said to be a \emph{graph epimorphism} if $f$ is surjective and given every edge $\{y_{1} , y_{2}\}$ of $Y$ there exists an edge $\{x_{1}, x_{2}\}$ of $X$ such that $f(x_{i}) = y_{i}$, $i \in \{1,2\}$.  A \emph{toroidal realization} of an abstract polyhedron $\P$ is a graph epimorphism $\beta: \sk(\P) \to \cX$ where $\cX$ is a graph in $\TT$. Following the convention used in Section\nobreakspace \ref {sec:basics}, any face of $\P$ can be identified with a subgraph of $\sk(\P)$. We say that a toroidal realization $\beta$ is \emph{faithful} if it is a graph isomorphism, and every face of $\P$ is uniquely determined by the corresponding subgraph of $\cX$. A toroidal realization is discrete if $V(\cX)$ is a discrete set in $\TT$. Due the compactness of $\TT$, a subset of points of $\TT$ is discrete if and only if it is finite. 

For the purposes of this work, we restrict our study to discrete toroidal realizations and unless otherwise specified, we will use the term \emph{toroidal realization} for short of \emph{discrete toroidal realization}. Given a toroidal realization of an abstract polyhedron $\P$, we define the \emph{group of toroidal symmetries} $G_{\bfL}(\P)$ as the group of isometries of $\TTL$ that induce automorphisms of $\P$. 

Note that every toroidal realization $\beta$ of an abstract polyhedron $\P$ induces a group homomorphism $\beta_{\ast}: \Delta \to G_{\bfL}(\P)$ of a certain subgroup $\Delta \leq \G (\P)$ in a similar way that an Euclidean realization of $\P$ induces an Euclidean representation of a certain subgroup of $\G(\P)$. We say that a toroidal realization is \emph{symmetric} if $\Delta= \G(\P)$. 

A natural problem is to classify symmetric and faithful toroidal realizations of regular polyhedra in $\TT$ to get analogous results to those presented in \cite{grunbaumOldAndNew,dress1,dress2,ordinary} for the Euclidean space $\E$ or those of \cite{roliProjective1,roliProjective2} for the projective space. 

Before attacking the problem presented above, we solve a slightly simpler problem: Consider a regular polyhedron $\P$ realized in the Euclidean space $\P$, that is, one of the 48 presented in Section\nobreakspace \ref {sec:basics}. For which lattice groups $\bfL$ the quotient map $\piL$ induces a symetric toroidal realization of a regular polyhedron $\PL$? We give an answer to this question in the following results but first we remark some general considerations and explain the general technique we use.

In the situation described above, we are interested in faithful discrete and symmetric toroidal realizations of regular polyhedra. Given a regular polyhedron $\P$ in the Euclidean space, we first classify those lattice groups $\bfL$ such that, in the case that $\P$ induces a discrete and faithful toroidal realization of a regular polyhedron $\PL$, then this realization will be symmetric. Latter we discuss on each case conditions for the symmetric toroidal realization to be discrete and faithful.

The general technique goes as follows: Suppose that a regular polyhedron $\P$ in $\E$ induces a discrete and faithful realization of a regular polyhedron $\PL$ in $\TTL$. If we want the realization to be symmetric we need that the isometries $R_0$, $R_1$ and $R_2$ that generate the symmetry group $G(\P)$ of $\P$ induce isometries in $\TTL$. According to Lemma\nobreakspace \ref {lem:equivNormGroup} this is equivalent to the extended special group $\bgop$ preserving the lattice $\L$. 

For every $i \in \{0,1,2\}$ let $R'_i$ denote the linear component of $R_i$. Now, according to Tables\nobreakspace \ref {tab:finite} and\nobreakspace  \ref {tab:pure}, for every finite regular polyhedron or regular pure apeirohedron in the Euclidean space $\P$, $R_0$, $R_1$ and $R_2$ are plane reflections or half-turns and so are $R'_0$ $R'_1$, $R'_2$. Then take $S_i \in \{R'_i , -R'_i\}$ such that $S_i$ is a plane reflection and consider de group $H(\P) := \left\langle S_0, S_1, S_2 \right\rangle \leq \bgop$. Since $G(\P)$ is a discrete group, then $\gop$ is a finite orthogonal group and so it is $\bgop$ (recall that $[\bgop : \gop] \leq 2$). Therefore $H(\P)$ is a finite irreducible group generated by reflections. 

In the case when $\P$ is a planar or a blended regular apeirohedron, $R_1$ and $R_2$ are two non-commuting involutions that fix de base vertex, therefore $R'_i$ for $i \in \{1,2\}$ is either a plane reflection or a half-turn. Considering again $S_i$ as before, the group $H(\P):= \left\langle S_1, S_2 \right\rangle$ is a reducible finite group generated by reflections. Finally observe that every lattice preserved by $\bgop$ must be preserved by $H(\P)$, then it is enough to know those lattices preserved by finite groups generated by reflections. 

We next study the possibilities for the group $H(\P)$. Since it is generated by reflections, if it is irreducible it must be isomorphic to $\Tet$, $\Oct$ or $\Ico$. If the group is reducible, then it is isomorphic to the dihedral group $D_n$ with $2n$ elements and it consists of the symmetries of a pyramid with base $\{n\}$. This group is generated by two reflections whose planes have angle $\frac{\pi}{n}$. We denote this group $\DDn$ since usually $D_n$ denotes the group generated by two half-turns at angle $\frac{\pi}{n}$. However, if we are interested in those groups that preserve lattices, the following result restricts the possibilities widely. A proof of the following theorem can be found in \cite[p. 152]{yaleGeomandSymm}, see also \cite[7E6]{ARP}.

\begin{thm}[Crystallographic Restriction]\label{thm:crystall}
 If $G$ is a group of isometries of $\bE^2$ or $\E$ that preserves a lattice, then $G$ does not contain rotations of periods other than $2$, $3$, $4$ or $6$.
\end{thm}

Theorem\nobreakspace \ref {thm:crystall} discards the groups $\Ico$ and $\DD_n$ for $n = 5$ and $n \geq 7$. Furthermore, since a lattice is preserved by an orthogonal group $G$ if and only if it is preserved by $\left\langle G, -id \right\rangle$, the lattices preserved by $\Tet$ are precisely those preserved by $\Oct$, since $\Oct = \langle \Tet, -id \rangle$. 

The following result gives a complete classification of the lattices preserved by $\Oct$. This result uses strongly the structure of lattices preserved by reflections given in Lemmas\nobreakspace \ref {lem:ReflectionLattices} and\nobreakspace  \ref {lem:bidimlattices}

\begin{lem}\label{lem:latticesOct}
 Let $\L$ be a lattice preserved by the group $\Oct$, then there exists $a \in \bR \sm \{0\}$ such that \[a^{-1} \L \in \{\cl, \bcl, \fcl \}.\]
 \begin{proof}
  Let $\L$ be as above. According to Table\nobreakspace \ref {tab:finite}, $R:=R_2$ and its conjugates $R'$ and $R''$ by $(R_1 R_0)$ and $(R_1 R_0)^{-1}$ respectively, are the reflections in the coordinated planes. Let $\Pi$ be the mirror plane of  $R_2$, that is, the plane of points whose third coordinate is zero.
  
  Let $\L_0 = \L \cap \Pi$. Since $\L$ is preserved by the reflection $R'$ on the plane $x=0$, $\L_0$ is preserved by the restriction of $R'$ to $\Pi$ and by Lemma\nobreakspace \ref {lem:bidimlattices} there exist a transformation $D = \diag(a,b)$ such that \[\L_0 D^{-1} \in \{\sql, \scl\}.\] Moreover, since $R_1 R_0$ permutes cyclically the coordinated planes then $a=b$ and thus $\L_0 D = a^{-1} \L_0$.  
  
  Suppose that $\L_0 = a \sql$. Then, by Lemma\nobreakspace \ref {lem:ReflectionLattices} there exists $w \in \L \sm \L_0$ such that $\L = \L_0 \+ w$ and $w$ projects to one among $o$, $a \frac{e_1}{2}$, $a \frac{e_2}{2}$ and $a\frac{e_1 + e_2}{2}$. However, since $R_0$ preserves the lattices the only possibilities are $o$ and $a\frac{e_1 + e_2}{2}$. Again, by the action of $R_1 R_0$, the former situation gives $\L = a \cl$ and the latter implies $\L = \frac{a}{2} \bcl$.
  
  If $\L_0 = a \scl$, take $w$ such that $w$ projects to one among $o$, $a \frac{e_1+e_2}{2}$, $a \frac{e_1-e_2}{2}$ and $a e_1$. Similar arguments that those used before show that unique possibility is that $w$ projects to $a e_1$ and this implies that $w= a(e_1 \pm e_3)$ and then $\L = a \fcl$. 
 \end{proof}
\end{lem}

It is straightforward to verify that the lattices $\cl$, $\bcl$ and $\fcl$ are actually preserved by $\Oct$. 

We proceed now to classify the lattices preserved by $\DD_2$, $\DD_3$, $\DD_4$ and $\DD_6$. According to the description of $\DDn$, this group is generated by two reflections $R, R'$ such that the angle between their mirrors is $\frac{\pi}{n}$. We may always assume, without loss of generality, that $R$ is the reflection in the plane $\Pi$ of points with first coordinate null and $R'$ is the reflection in the plane with the appropriate angle and such that $\Pi \cap \Pi'$ is the $z$ axis. With those assumptions it is clear that $\DD_2 \leq \DD_4$ and $\DD_3 \leq \DD_6$, therefore it is enough to find those lattices preserved by $\DD_2$ and $\DD_3$ and determine which of them are preserved by $\DD_4$ and $\DD_6$ respectively.

The following lemma classifies the lattices preserved by $\DD_{2}$.

\begin{lem}\label{lem:latticesD2}
 Let $\L$ be a lattice preserved by $\DD_2$ then exists a diagonal transformation $D$ such that \[\L D^{-1} \in \{\cl, \bcl, \fcl, \scl \+ e_3, 2\sql \+ (e_1+e_3), 2\sql \+ (e_2 + e_3)\}.\]
  \begin{proof}
    Let $R$ and $R'$ be the reflections described above. Let $\Sigma$ be the plane $z=0$. Since $\L$ is preserved by $R$ and $R'$, then it is also preserved by the half-turn $R R'$.  By Lemma\nobreakspace \ref {lem:RefequivRot}, $\L$ is preserved by the reflection through $\Sigma$. Let $\L_0 = \L \cap \Sigma$. Since $\L_0$ is preserved by the Lemma\nobreakspace \ref {lem:bidimlattices} there exists a diagonal transformation $D_0 = \diag(d_1,d_2)$ such that $\L_0 D_{0}^{-1} \in \{ \sql, \scl\}$

    If $\L_0 D^{-1} = \sql$ then $\{d_1 e_1, d_2 e_2 \}$ is a basis of $\L$ and we can take $w$ such that $w$ projects to exactly one among $o$, $\frac{d_{1}}{2} e_{1}$, $\frac{d_{2}}{2} e_{2}$ or $\frac{d_1}{2} e_{1} + \frac{d_{2}}{2} e_{2}$ and $\L = \L_0 \+ w$. If $w$ projects to $o$, then $w = d_3 e_3$ and then $\L = \cl D$ with $D= \diag(d_{1}, d_{2}, d_{3})$. If $w$ projects to $\frac{d_{1}}{2} e_{1}$ then $w =  \frac{d_{1}}{2} e_{1} + \frac{d_{3}}{2} e_{3}$ for some suitable $d_{3}$; define $D =\diag(\frac{d_{1}}{2},\frac{d_{2}}{2},\frac{d_{3}}{2})$ and then $\L = (2\sql \+ (e_{1}+ e_{3}))D$. Similarly we can conclude that if $w$ projects to $\frac{d_{2}}{2} e_{2}$, then $\L = (2\sql \+ (e_{2}+ e_{3}))D$. Finally, if $w$ projects to $\frac{d_1}{2} e_{1} + \frac{d_{2}}{2} e_{2}$, then $w=\frac{d_1}{2} e_{1} + \frac{d_{2}}{2} e_{2} + \frac{d_{3}}{2} e_{3}$ for some $d_3 \in \bR$, in which case define $D =\diag(\frac{d_{1}}{2},\frac{d_{2}}{2},\frac{d_{3}}{2})$ and then $\L = \bcl D^{-1}$.
    
    Similar arguments prove that if $\L_0 =  \L_{(1,1)} D_0$, then $\L D$ is either $\fcl$ or $\L_{(1,1)} \+ e_3$.
  \end{proof}
\end{lem}

It is easy to verify that all the lattices mentioned in Lemma\nobreakspace \ref {lem:latticesD2} are preserved by $\DD_2$. As we said before, we can use the lemma above to classify the lattices preserved by $\DD_4$, since every lattice preserved by $\DD_4$ must be preserved by $\DD_2$. First observe that there is an element of $\DD_4$ that swaps $e_1$ and $e_2$, this discards the lattices $2\sql \+ (e_1+e_3)$ and $2\sql \+ (e_2 + e_3)$. Moreover, this forces that the transformation $D$ mentioned in Lemma\nobreakspace \ref {lem:latticesD2} has the form $D= \diag(a,a,b)$ for some $a,b \in \bR$. We summarize these observations in the following result.

\begin{lem}\label{lem:latticesD4}
If $\L$ is a lattices preserved by the group $\DD_4$ then there exists a diagonal transformation $D=(a,a,b)$ for some $a,b > 0$ such that \[\L D^{-1} \in \left\{\cl, \bcl, \fcl, \scl \+ e_3\right\}. \]
\end{lem}

It is not hard to verify that all those lattices are preserved by $\DD_4$. With the notation of Lemma\nobreakspace \ref {lem:latticesD4}, the lattice $\cl D$ is the lattice determined by the vertices of a tessellation with prisms of height $b$ over squares of side $a$ parallel to the coordinated axes. The lattices $\bcl D$ and $\fcl D$ consist of the vertices of a tessellation with prisms with height $2b$ and squared base of side $2a$ together with the centres of the prisms in the former case and the centres of the faces of the prisms in the latter. The lattice $(\scl \+ e_3)D$ consists of the vertices of a tessellation with prisms of height $b$ over squares of side $\sqrt{2} a$ whose edges make an angle of $\frac{\pi}{4}$ with the coordinated axes.

The analogous classification result for the group $\DD_3$ is a slightly more complicated, so we have decided to split it in some small lemmas. At this moment it is convenient to denote $u_{1} = e_1$ and $u_2 = \frac{1}{2}e_1 + \frac{\sqrt{3}}{2} e_2$ such that $\TB=\{u_1, u_{2}\}$ is a basis for the lattice $\tl$.

\begin{lem}\label{lem:latticesD3bidim}
 Let $\L$ be a rank-$2$ lattice in $\bE^2$, $R$ the reflection with respect to the line $\Pi$ that contains $o$ and $u_1$ and $R'$ the reflection in the line $\Pi'$ that contains $o$ and $u_2$. If $\L$ is preserved by $R$ and $R'$, then there exists $a > 0$ such that \[a \L \in \left\{\tl, \tcl\right\}.\]
 \begin{proof}
  Let $v_1 \in \Pi \cap \L$ and $v_2 \in \Pi' \cap \L$ with minimum positive distance to $o$ such that the angle between them is $\frac{\pi}{3}$. Observe that by minimality of $|v_1|$ and $|v_2|$, $v_2 +v_2 R = v_1$. Therefore the triangle $o v_1 v_2$ is isosceles and has an angle of $\frac{\pi}{3}$ and hence it is equilateral and $|v_1|=|v_{2}|$. If $\{v_1 , v_2\}$ is a basis for $\L$, then $a\L = \tl$ for some $a>0$.
  
  If there is a point of $\L$ that is not an integer linear combination of $v_1$ and $v_{2}$, then there is a point $x \in \L$ in the parallelogram spanned by $v_1$ and $v_2$. Moreover, we may assume that $x$ is in the triangle $o v_{1} v_{2}$, because if $x$ is in the triangle determined by $v_1 $, $v_2$ and $ v_1+v_2$, then $xR + v_{2} - v_{1} $ belongs to the triangle $o v_{1} v_{2}$. Observe that $x + xR \in \Pi$ but, by the choise of $x$ and the minimality of $|v_{1}|$ $x + xR = v_1$ and then $x$ projects to $\Pi$ to the midpoint of $o$ and $v_1$; analogously we may conclude that $x$ projects to $\Pi'$ to the midpoint of $o$ and $v_2$. Therefore, $x$ must be the centre of the triangle $o v_{1} v_{2}$ and then $a\L = \tcl$ for some $a > 0$.
 \end{proof}
\end{lem}

\begin{lem}\label{lem:latticesD3sublat}
 Let $\L$ be a rank-$3$ lattice preserved by a $3$-fold rotation $S$ on a line $L$ that contains the origin $o$. Let $\Pi$ the orthogonal plane to $\L$ that contains $o$. Then $\L \cap L$ is a rank-$1$ lattice and $\L \cap \Pi$ is a rank-$2$ lattice. Furthermore, if $x \in \L$, then $3x$ projects to $\Pi$ in a point of $\L \cap \Pi$.
 \begin{proof}
  Observe that if $x \in \L \sm \{o\}$ then $x + xS + xS^2 \in \L \cap L$, then take a point $y$ closest to $o$ among all the points of $(\L \cap L) \sm \{o\}$. Observe that0 $y$ generates $\L \cap L$. Similarly, if $x \in \L \sm L $, then $x - xS$ and $x - xS^{2}$ are two linearly independent vectors in $\Pi \cap \L$; among all pairs of linearly independent vectors of $\Pi \cap \L$ take $\{y_1, y_2\}$ such that the parallelogram spanned by $y_{1}$ and $y_{2}$ has minimum area and then $\{y_{1}, y_{2}\}$ is a basis for $\L \cap \Pi$.
  
  For any $x \in \L$, $3x - (x + xS + xS^2)$ is the projection of $3x$ to $\Pi$, which proves the last part of the lemma.
 \end{proof}
\end{lem}

Given the Lemmas\nobreakspace \ref {lem:latticesD3bidim} and\nobreakspace  \ref {lem:latticesD3sublat} it is now easy to prove the following result, which classifies the lattices preserved by the group $\DD_3$.

\begin{lem}\label{lem:latticesD3}
 Let $\L$ be a rank-$3$ lattice. If $\L$ is preserved by the group $\DD_3$, then there exists a diagonal transformation $D= \diag(a,a,b)$ such that \[ \L D^{-1} \in \left\{ \tl \+ e_3, \tcl \+ e_3 , \tcl \+ (u_1 + e_3), \tcl \+ (u_2 + e_3)\right\}\]
 \begin{proof}
  Assume that $\L$ is as above and that the group $\DD_3$ is generated by the reflection $R$ in the plane $\Pi$ that contains $o$, $u_1$ and $e_3$ and the reflection $R'$ in the plane $\Pi'$ that spanned by $o$, $u_2$ and $e_3$. Let $\Sigma$ be the plane $z=0$ and $L = \Sigma^{\perp}$, that is, the $z$-axis. Observe that $RR'$ is a $3$-fold roatation in $L$ and by Lemma\nobreakspace \ref {lem:latticesD3sublat}, $\L_0 :=\L \cap \Sigma$ is a rank-$2$ lattice. Observe that $\L_0$ is preserved by the restriction of $R$ and $R'$ to the plane $\Sigma$, hence by Lemma\nobreakspace \ref {lem:latticesD3bidim} there exists $a >0$ such that $ a^{-1} \L  \in \left\{ \tl, \tcl \right\}$.
  
  Assume that $ \L  = a\tl$. By Lemma\nobreakspace \ref {lem:latticesD3sublat}, $\L \cap L$ is a lattice of rank $1$. Let $b >0$ such that $be_3$ is the closest point of $\L\cap L$ to $o$. If $\L = \L_0 \+ be_3$ define $D= \diag(a,a,b)$ and then $\L D^{-1} = \tl \+ e_3$. Observe that for any point $x \in \Pi^{\perp}$, $d(x, \Pi) \geq \sqrt{3}a$. Suppose that $\L \sm (\L_0 \+ be_3) \neq \emptyset$. In this situation there must exist a point $x \in \L$ that projects to $\Sigma$ in the triangle determined by $o$,$ au_1 $ and $ au_2$, such that $d(x, \Sigma) < b$. If $x \not\in \Pi$, then $x-xR \in \Pi^{\perp}$ and $|x-xR| < \sqrt{3}a$ which is a contradiction, therefore $x \in \Pi$. We can proceed in the analogous way to conclude that $x \in \Pi'$, then $x \in \L \cap L$ which contradicts the minimality of $b$. Therefore $\L = \L_0 \+ be_3$. 
  
  Suppose then that $\L_0 = a \tcl$ and take $be_3$ as before. If $\L = \L_0 \+ be_3$ then $\L D^{-1} = \tcl \+ e_3$ for $D= \diag(a,a,b)$. Otherwise, there must exist some points of $\L$ that project to $\Sigma$ in the parallelogram $P$ determined by $o$, $au_1 +au_2$ and $2au_1-au_2$. Among all of them, take $x$ closest to $\Sigma$ such that $\L = \L_0 \+ x$. Since $x$ projects to $P$, then $d(x, \Pi) \leq \frac{\sqrt{3}}{2}a$. Recall that $x-xR \in \Pi^{\perp}$ so either $x \in \Pi$ or $d(x,\Pi)=\frac{\sqrt{3}}{2}a$ but the latter implies that one of $x$ or $xR$ projects to $au_1 + au_{2}$ which contradict the assumption of $\L$ not being a vertical translation lattice, therefore $x \in \Pi$. Recall that by Lemma\nobreakspace \ref {lem:latticesD3sublat}, $3x$ projects to a point of $\L_0$, so $x$ must project either to $au_1$ or $2a u_1$ and by the choice of $x$, $d(x,\Sigma) \in \left\{ \frac{b}{3}. \frac{2b}{3}\right\}$. Suppose then that $x$ projects to $au_1$, if $d(x, \Sigma)=\frac{2b}{3}$ then $d(2x-be_3, \Sigma) = \frac{b}{3}$, which contradicts the choice of $x$. Similarly, if $x$ projects to $2a u_1$ and $d(x, \Sigma)=\frac{2b}{3}$ then $d(2x-3au_1 -be_3, \Sigma) = \frac{b}{3}$. Therefore, in any case $d(x, \Sigma)$ must be $\frac{b}{3}$. Finally observe that $\L_0 \+ (2au_1 + \frac{b}{3}e_3) = \L_0 \+ (au_2 + \frac{b}{3}e_3)$ and thus we have that \[\L D^{-1} \in \left\{ \tcl \+ (u_1 + e_3), \tcl \+ (u_2 + e_3) \right\}\] for $D = \diag(a,a,\frac{b}{3})$ as desired.
 \end{proof}
\end{lem}

It is straightforward to verify that all lattices listed in Lemma\nobreakspace \ref {lem:latticesD3} are preserved by $\DD_3$. This result together with the fact that $\DD_3 \leq \DD_6$ gives us an easy prove of the classification result for those lattices preserved by $\DD_6$. 

\begin{lem}\label{lem:latticesD6}
 Let $\L$ be a lattice. If $\L$ is preserved by $\DD_6$ then there exists a diagonal transformation $D =\diag(a,a,b)$ such that \[\L D^{-1} \in \left\{ \tl \+ e_3, \tcl \+ e_3 \right\}.\]
 \begin{proof}
  Any lattice preserved by $\DD_6$ must be preserved by $\DD_3$. The result follows form the fact that $\DD_6$ does not preserve the lattices $\tcl \+ (u_1 + e_3)$ and $\tcl \+ (u_2 + e_3)$ but it does preserve the lattices $\tl \+ e_3$ and $ \tcl \+ e_3$.
 \end{proof}
\end{lem}

According to the discussion at the beginning of the section, we are interested in determine those lattice groups $\bfL$ such that a regular polyhedron $\P$ in $\E$ induces a regular polyhedron $\PL$ in $\TTL$. We proved that a necessary condition is that the lattice $\L$ is preserved by the extended special group $\bgop$. To determine those lattices preserved by $\bgop$ we associate to each polyhedron a group $H(\P)\leq \bgop$ generated by reflections, study the possibilities for $H(\P)$ and determine all the lattices preserved by each of these possibilities.

In Table\nobreakspace \ref {tab:groupH} we list the 48 regular polyhedra in $\E$ according to the group $H(\P)$ described before.

\begin{table}
\begin{center}
\resizebox{\textwidth}{!} {
{%
\newcommand{\mc}[3]{\multicolumn{#1}{#2}{#3}}
\begin{tabular}{|c|cc|cc|c|c|cc|}
\hline
$\Tet$ & \mc{2}{c|}{$\Oct$} & \mc{2}{c|}{$\Ico$} & $\DD_3$ & $\DD_4$ & \mc{2}{c|}{$\DD_6$}\\ \hline
$\{3,3\}$ & $\{3,4\}$ & $\{4,3\}$ & $\{3,5\}$ & $\{5,3\}$ & $\{3,6\}\# \{\infty\}$ & $\{4,4\}$ & $\{3,6\}$ & $\{6,3\}$\\
$\{6,3\}_{4}$ &$\{4,3\}_{3}$ & $\{6,4\}_{3}$ &  $\{10,5\}$  & $\{10,3\}$ & $\{\infty,6\}_{3} \# \{\infty\}$ & $\{\infty,4\}_{4}$ & $\{\infty,6\}_{3}$ & $\{infty,3\}_{6}$\\
$\{6,6|3\}$ & $\{6,4|4\}$ & $\{4,6|4\}$ & $\{6,\frac{2}{2}\}$ & $\{5,\frac{5}{2}\}$ & $\{6,3\} \# \{\ \} $ & $\{4,4\} \# \{\ \}$ & $\{3,6\} \# \{\ \}$ & $\{\infty,6\}_{3} \# \{ \}$ \\
$\{\infty,6\}_{4,4}$ & $\{\infty,4\}_{6,4}$ & $\{\infty,6\}_{6,3}$ & $\{\frac{5}{2},5\}$ & $\{6,5\}$ & $\{\infty,3\}_{6} \# \{\ \}$  & $\{\infty,4\}_{4} \# \{ \}$ & $\{6,3\} \# \{\infty \} $ & $\{\infty,3\}_{6} \# \{\infty\}$\\
$\{6,6\}_{4}$ & $\{6,4\}_{6}$ & $\{4,6\}_{6}$ & $\{\frac{10}{3},3\}$ & $\{\frac{5}{2},3\}$ &  & $\{4,4\} \# \{\infty\}$ &  & \\
$\{\infty,3\}^{(a)}$ & $\{\infty,4\}_{\cdot, \ast 3}$ & $\{\infty,3\}^{(b)}$ & $\{3,\frac{5}{2}\}$ & $\{\frac{10}{3},\frac{5}{2}\}$ &  & $\{\infty,4\}_{4} \# \{ \}$ &  & \\ \hline
\end{tabular}
}%

}
\caption{Regular polyhedra in $\E$ according to its group $H(\P)$.} \label{tab:groupH}
\end{center}
\end{table}

It is important to remark that in the case of the polyhedra $\{6,3\} \# \{\ \}$ and its Petrial $\{\infty,3\}_{6} \# \{\ \}$ the group $H(\P)$ is not  actually $\DD_3$ but its conjugate by a rotation of $\frac{\pi}{2}$ such that the generating reflections mentioned before are the reflections on the plane $x=0$ and in the plane $x+\sqrt{3y} = 0$. Hence the latices preserved by this group are the result of rotating $\frac{\pi}{3}$ those lattices preserved by $\DD_3$.

Now we can give some results about the classification of regular polyhedra in $\TT$. Theorem\nobreakspace \ref {thm:Ico} is a direct consequence of Theorem\nobreakspace \ref {thm:crystall}, since $\Ico$ contains a $5$-fold rotation.

\begin{thm}\label{thm:Ico}
 Let $\P$ be a finite regular polyhedra in $\E$ whose symmetry group is $\Ico$. Then there are no lattice groups $\bfL$ such that the quotient map $\piL$ induces a symmetric toroidal realization of a regular polyhedron. 
\end{thm}

Theorem\nobreakspace \ref {thm:TetandOct} gives necessary conditions to the existence of regular polyhedra in $\TTL$ induced by regular polyhedra in $\E$ whose symmetry group is $\Tet$ or $\Oct$. This result follows from Lemma\nobreakspace \ref {lem:latticesOct}.

\begin{thm}\label{thm:TetandOct}
 Let $\P$ be a regular polyhedron in $\E$ such that $G(\P)$ is either $\Tet$ or $\Oct$ and let $\bfL$ be a lattice group. If the quotient map $\piL$ induces and symmetric toroidal realization of a regular polyhedron in $\TTL$, then there exists $a >0$ such that \[a^{-1} \bfL \in \left\{ \Cl, \Fcl, \Bcl \right\}.\] 
\end{thm}

If we are interested in regular polyhedra in $\TTL $ induced by a regular apeirohedron in $\E$, Lemmas\nobreakspace  \ref {lem:latticesOct} to\nobreakspace  \ref {lem:latticesD4} ,  \ref {lem:latticesD3} and\nobreakspace  \ref {lem:latticesD6} and the information in Table\nobreakspace \ref {tab:groupH} impose necessary conditions over the lattice group $\bfL$. This conditions are summarized in the following theorem.

\begin{thm}\label{thm:infinite}
	Let $\P$ be a regular apeirohedron in $\E$. Suppose that $\bfL$ is a lattice group such that the quotient map $\piL$ induces a symmetric toroidal realization of a finite regular polyhedron $\PL$, then we have the following possibilities:
	\begin{enumerate}[label=\textit{(\roman*)}]
	 \item \label{inc:pure} If $\P$ is a pure apeirohedron, then there exists $a >0$ such that \[a^{-1} \bfL \in \left\{ \Cl, \Fcl, \Bcl \right\}.\]
	 \item \label{inc:D4} If $\P$ is one among $\{4,4\}$, $\{4,4\} \# \{\ \}$, $\{4,4\} \# \{\infty \}$ or the Petrial of one of these, then there exists a diagonal transformation $D=\diag(a,a,b)$ such that \[\bfL D^{-1} \in \left\{\Cl, \Fcl, \Bcl, \left\langle \Scl, t_3 \right\rangle \right\}. \] 
	 \item \label{inc:D3a} If $\P$ is the polyhedron $\{3,6\} \# \{\infty\}$ or its Petrial, the polyhedron $\{\infty,6\}_{3} \# \{\infty\}$ there exists a diagonal transformation $D= \diag(a,a,b)$ such that \[ \bfL D^{-1} \in \left\{ \left\langle \Tl, t_3 \right\rangle, \left\langle \Tcl, t_3 \right\rangle , \left\langle \Tcl, t_1 \right\rangle, \left\langle \Tcl, t_{2}\right\rangle \right\},\] where $t_1$ is the translation by $(u_1 + e_3)$, $t_2$ is the translation by $(u_2 + e_3)$ and $t_3$ is the translation by $e_3$, with $\{u_1, u_2\}$ the basis of $\tl$.
	 \item \label{inc:D3b} If $\P$ be the polyhedron $\{6,3\} \# \{\ \}$ or its Petrial, the polyhedron $\{\infty,3\}_{6} \# \{\ \}$, then there exists a diagonal transformation $D= \diag(a,a,b)$ such that \[ \bfL D^{-1} \in \left\{ \left\langle \Tl, t_3 \right\rangle, \left\langle \Tcl, t_3 \right\rangle , \left\langle 3\Tl, t_1 \right\rangle, \left\langle 3\Tl, t_{2}\right\rangle \right\},\] where $t_1$ is the translation by $u_1 + u_2 + e_3$, $t_2$ is the translation by $2u_1-u_2 + e_3$ and $t_3$ is the translation by $e_3$, with $\{u_1, u_2\}$ the basis of $\tl$. 
	 \item \label{inc:D6} If $\P$ is one among $\{3,6\}$, $\{6,3\}$, $\{3,6\} \# \{\ \}$, $\{6,3\} \# \{\infty \}$ or the Petrial of one of these, then there exists a diagonal transformation $D= \diag(a,a,b)$ such that \[\bfL D^{-1} \in \left\{\left\langle \Tl, t_3 \right\rangle,\left\langle \Tcl, t_3 \right\rangle \right\}, \] where $t_3$ is the translation by $e_3$.
	\end{enumerate}
\end{thm}

It is important to emphasize that the sufficiency of such results depend on the choice of one (in the case of finite polyhedra or pure apeirohedra) or two (in the case of planar or blended apeirohedra) parameters. The election of those parameters will determine the discreteness of the vertex-set of $\PL$ and whether or not the diamond condition holds. We will discuss this in the following subsection.

\subsection{Determining the parameters.}

The aim of this subsection is to determine possible values for the parameters in Theorems\nobreakspace \ref {thm:TetandOct} and\nobreakspace  \ref {thm:infinite} such that the corresponding quotient map induces a toroidal realization of a regular polyhedron $\PL$. 
First observe that the lattices listed in such results are all preserved by the corresponding extended special group. So, given appropriate parameters that guarantee a toroidal realization of any polyhedron, such realization will be symmetric.

\subsubsection{Finite polyhedra}
According to Theorems\nobreakspace \ref {thm:Ico} and\nobreakspace  \ref {thm:TetandOct} we only need to verify when the quotient map $\piL$ induces a toroidal realization of the tetrahedron, the cube, the octahedron and their Petrials when $\bfL$ is one of the lattices listed in Theorem\nobreakspace \ref {thm:TetandOct}.

In the situation described above, there is one parameter $a$ to determine. We are interested in faithful and symmetric realizations, but as mentioned before, symmetry is given by restricting the possibilities of $\bfL$ to scalar multiples of $\Cl$, $\Fcl$ and $\Bcl$. If we want the realization to be faithful, a necessary condition is that no two vertices of the polyhedron $\P$ are equivalent under the action of $\bfL$. Furthermore, according to the definition of toroidal realization we do not want that a vertex of $\P$ identifies with an interior point of any edge of $\P$. This offers a trivial lower bound $\alpha_{0}$ for $a$ given by the property that if $a \leq \alpha_{0}$ then there is a vertex such that it is equivalent to a point of some incident edge.

There is also another bound $\alpha_{1}$ for $a$ given by the property that if $\alpha_{1} < a$ then the $\P$ is properly contained in $D(\bfL)$, the Dirichlet domain of $\bfL$ and therefore no two points of $\P$ are equivalent under the action of $\bfL$. 

It is possible to determine which of the values of $a$ in $(\alpha_{0},\alpha_{1}]$ give faithful toroidal realizations of $\P$ by few computations. Observe that given that $\alpha_{0} < a$ there are a finite (rather small) number of translations of $\bfL$ such that their translation vector have length smaller that the diameter of $\P$ (bear in mind that $\P$ is a compact set of $\E$). Therefore, we just have to verify that those translations do not give undesirable identifications. This can be computed by hand in just one vertex of $\P$ since $G(\P)$ acts transitively on vertices and preserves $\bfL$. 

Just to show an example, if $\P$ is the octahedron with edge length $\sqrt{2}$ and $\bfL= \Cl$ then $\alpha_{0} = 1$ and $\alpha_{1} = 2$. If $\alpha_{0} < a \leq  \alpha_{1}$ the only translations that have to be verified are those given by the vectors $(\pm a,0,0)$, $(0,\pm a,0)$, $(0,0,\pm a)$, $(\pm a, \pm a,0)$, $(0,\pm a,\pm a)$ and $(\pm a,0,\pm a)$, since any other translation vector in $\L$ will have length greater than $2$, which is the diameter of $\P$. It is easy to verify that those translations mentioned before will not give place to undesired identifications unless $a=\alpha_{1}=2$.

In Table\nobreakspace \ref {tab:parametersFinite} we show the values of $\alpha_{0}$, $\alpha_{1}$ and the valid values of $a$ in $(\alpha_{0}, \alpha_{1}]$ for the tetrahedron, the cube and the octahedron and each of the lattices preserved by their symmetry group. Their Petrials are no listed since they share the $1$-skeleton and therefore, the values are the same.

Theorems\nobreakspace \ref {thm:Ico} and\nobreakspace  \ref {thm:TetandOct} and Table\nobreakspace \ref {tab:parametersFinite} complete the classification of regular polyhedra in $\TT$ induced by finite polyhedra in $\E$. It is important to remark that we do not impose the condition of the realizations being faithful, but weaker conditions force all of them to be faithful anyway.

\begin{table}
\begin{center}
{%
\newcommand{\mc}[3]{\multicolumn{#1}{#2}{#3}}
\begin{tabular}{|c|c|c|c|c|c|c|c|c|c|}
\hline
$\P$ & \mc{3}{c|}{$a\Cl$} & \mc{3}{c|}{$a \Fcl$} & \mc{3}{c|}{$a\Bcl$}\\
 & \mc{1}{c}{$\alpha_{0}$} & \mc{1}{c}{$\alpha_{1}$} & \mc{1}{c|}{$(\alpha_{0}, \alpha_{1}]$} & \mc{1}{c}{$\alpha_{0}$} & \mc{1}{c}{$\alpha_{1}$} & \mc{1}{c|}{$(\alpha_{0}, \alpha_{1}]$}& \mc{1}{c}{$\alpha_{0}$} & \mc{1}{c}{$\alpha_{1}$} & \mc{1}{c|}{$(\alpha_{0}, \alpha_{1}]$}\\ \hline
$\{3,3\}$ & $2$ & $2$ & -- & $2$ & $2$ & -- & $1$ & $2$ & $(1,2]$\\
$\{3,4\}$ & $1$ & $2$ & $(1,2)$ & $1$ & $1$ & -- & $\frac{1}{2}$ & $1$ & $\left(\frac{1}{2},1\right)$ \\ 
$\{4,3\}$ & $2$ & $2$ & -- & $1$ & $2$ & $(1,2)$ & $1$ & $2$ & $(1,2)$ \\\hline
\end{tabular}
}%
\caption{Parameters for finite polyhedra} \label{tab:parametersFinite}
\end{center}
\end{table}

\subsubsection{Pure apeirohedra.}

In this subsection we determine the parameter $a$ for each of the lattices mentioned in part \ref{inc:pure} of Theorem\nobreakspace \ref {thm:infinite} for each pure apeirohedron $\P$. As mentioned before it is not possible to have discrete an faithful toroidal realizations of regular apeirohedra, however we may ask for those parameters such that the quotient map $\piL$ induce faithful and discrete toroidal realization of a finite polyhedron $\PL$.

If we are interested in toroidal realizations we cannot admit a vertex of $\P$ being identified with an interior point of an edge. If we want the vertex figures of $\PL$ to be polygons, if two vertices of $\P$ are identified then they must be equivalent under the action of the group $\T(\P)$ of translations of $\P$. Furthermore, if we want the faces of $\PL$ to be polygons we need that if two vertices in the same face are identified, then the translation that identifies them must map the face to itself.

The previous observations gives us a technique to determine the parameter $a$ for eight of the twelve polyhedra, namely all but $\{6,4 | 4\}$, $\{\infty,3 \}^{(b)}$ and their Petrials. If $\P$ is not $\{6,4 | 4\}$, $\{\infty,3 \}^{(b)}$ or  the Petrial of one of these, then there exists a line $L$ in $\E$ such that every point of $L$ belongs to an edge of $\P$ that is contained in $L$. If $v$ is a vertex of $\P$ in $L$ then every translation of $\bfL$ in the direction of $L$ must identify $v$ with other vertex $w$ in $L$ such that $v-w$ is the translation vector of an element in $\T(\P)$. 

The restriction mentioned above together with the restriction that no two non-translate vertices of $\P$ will be identified are enough to determine the possible values of $a$ for such polyhedra. Small values of $a$ must be ruled out in some cases since they identify two vertices in a finite face, which implies that faces of $\PL$ are not cycles. The values obtained are listed in Table\nobreakspace \ref {tab:parametersPure}.

If $\P$ is $\{6,4 | 4\}$, $\{\infty,3 \}^{(b)}$ or  the Petrial of any of these a more detailed analysis is required. If we want $\PL$ to be finite, then $a$ must be a rational number $\frac{p}{q}$ with $p,q \in \bN$ such that $\gcd(p,q) = 1$. We explain the technique with $\P = \{\infty,3 \}^{(b)}$ and $\bfL = a \Fcl$, similar arguments will be useful to determine the parameters for the other two lattices and $\P$ being one among $\{\infty,3 \}^{(a)}$, $\{6,4 | 4\}$ and $\{\infty, 4\}_{6,4}$.

The base edge of $\P$ is the segment determined by $o$ and $(1,1,0)$. Since the translation by $(4,4,0)$ belongs to $G(\P)$, if there is a multiple of $\frac{p}{q}$ of the form $4n + \epsilon$ with $n \in  \bN$ and $\epsilon \in (0,1]$, then the corresponding translation will identify the base vertex $o$ with an interior point of an edge or with a vertex which is not a translate of $o$ by an element of $\T(\P)$. This rules out any value of $\frac{p}{q}$ such that $q \geq 4$, since there are integers $x,y$ such that $xp+yq=1$ which implies that $4x \frac{p}{q} = -4y + \frac{4}{q}$.

We only have to explore the cases $q \in \{1,2,3\}$. If $\gcd(p,4q) =1$ there exist $x,y \in \bZ$ such that $xp +4yq = 1$ which implies $x \frac{p}{q} = -4y + \frac{1}{q}$. It only rests to verify $a = p$ with $p$ even and $a= \frac{p}{3}$ with $p \equiv \pm 2, \pm 4 \pmod{12}$, however $a=\frac{12k \pm 2}{3}$ implies $a = 4k \pm \frac{2}{3}$. If $a=2k$  then $a\Fcl \leq \T(\P)$ and it is easy to verify that for any $k$ the quotient map induce a toroidal realization. Observe that if $a= \frac{4}{3}$ then no undesired identifications occur, since the lattice $\frac{4}{3}\fcl$ only intersect $\P$ in vertices with coordinates $(4m,4n,4k)$, which are equivalent under $\T(\P)$; this implies that for $a= \frac{12k\pm4}{3}$ with $k \in \bN \cup \{0\}$, the quotient map induces toroidal realizations of a finite polyhedron $\{\infty,3 \}^{(b)}_{a\Fcl}$ in $a\Fcl$, since $a \Fcl \leq \frac{4}{3} \Fcl$.

\begin{table}
\begin{center}
\resizebox{\textwidth}{!} {

\begin{tabular}{|c | c | c | c| c|}
	\hline 
	$\P$ & $\T(\P)$ & \multicolumn{3}{c|}{Possible values of $a$} \\ 
	&  & \multicolumn{1}{c}{$\bfL = a \Cl$} & \multicolumn{1}{c}{$\bfL = a \Fcl$} & \multicolumn{1}{c|}{$\bfL = a \Bcl$} \\ 
	\hline
	$\{4,6|4\}$ & $\Bcl$ & $a \in 2\bN$ & $a \in 2\bN$  & $a \in \bN$ \\ \hline
	$\{6,4|4\}$ 	& $\Bcl$ 	& $a \in 2\bN$, 	& $a \in 2\bN $, 	& $a \in \bN, $ \\
	&	& $a= \frac{p}{3}, \ p\equiv \pm 4 \pmod{12} $&	$a= \frac{p}{3}, \ p\equiv \pm 4 \pmod{12} $& $a= \frac{p}{6}, \ p\equiv \pm 4 \pmod{12} $\\\hline
	$\{6,6|3\}$ & $2\Fcl$ & $a \in 2(\bN \sm \{1\})$ & $a \in 2(\bN \sm \{1\})$ & $a \in \bN \sm \{1\}$ \\ \hline
	$\{\infty,6\}_{4,4}$ & $\Bcl$ & $a \in 2\bN$ & $a \in 2\bN$ & $a \in \bN$ \\ \hline
	$\{ \infty,4\}_{6,4}$ 	& $\Bcl$ 	& $a \in 2\bN$, 	& $a \in 2\bN $, 	& $a \in \bN, $ \\
	&	& $a= \frac{p}{3}, \ p\equiv \pm 4 \pmod{12} $&	$a= \frac{p}{3}, \ p\equiv \pm 4 \pmod{12} $& $a= \frac{p}{6}, \ p\equiv \pm 4 \pmod{12} $\\\hline
	$\{\infty,6\}_{6,3}$ & $2\Fcl$ & $a \in 2\bN$ & $a \in 2\bN$ & $a \in \bN$ \\ \hline
	$\{6,6\}_{4}$ & $2 \Cl$ & $a \in 2(\bN\sm\{1\})$ & $a \in 2\bN$ & $a \in \bN$ \\ \hline
	$\{4,6\}_{6}$ & $2 \Cl$ & $a \in 2\bN$ & $a \in 2\bN$ & $a \in \bN$ \\ \hline
	$\{\infty,3\}^{(b)}$ 	& $2\Bcl$ 	& $a \in 2\bN$, 	& $a \in 2\bN $, 	& $a \in \bN, $ \\
	&	& $a= \frac{p}{3}, \ p\equiv \pm 4 \pmod{12} $&	$a= \frac{p}{3}, \ p\equiv \pm 4 \pmod{12} $& $a= \frac{p}{6}, \ p\equiv \pm 4 \pmod{12} $\\\hline
	$\{\infty,3\}^{(a)}$ &  $2\Bcl$ 	& $a \in 2\bN$, 	& $a \in 2\bN $, 	& $a \in \bN, $ \\
	&	& $a= \frac{p}{3}, \ p\equiv \pm 4 \pmod{12} $&	$a= \frac{p}{3}, \ p\equiv \pm 4 \pmod{12} $& $a= \frac{p}{6}, \ p\equiv \pm 4 \pmod{12} $\\\hline
	$\{6,4\}_{6}$ & $2\Cl$ & $a \in 2\bN$ & $a \in \bN \sm \{1\}$ & $a \in \bN\sm \{1\}$ \\ \hline
	$\{\infty, 4\}_{\cdot, \ast 3}$ & $\Bcl$ & $a \in 2\bN$ & $a \in \bN \sm \{1\}$ & $a \in \bN$ \\ 
	\hline
	
\end{tabular} 

}
\caption{Parameters for pure apeirohedra} \label{tab:parametersPure}
\end{center}
\end{table}

\subsubsection{Planar and blended apeirohedra.}

According to Theorem\nobreakspace \ref {thm:infinite} if $\P$ is a planar or blended regular apeirohedra such that $\P$ induces a symmetric toroidal realization of a regular polyhedron $\PL$, there are two parameters $a$ and $b$ to determine. The lattices preserved by the groups of planar and blended apeirohedra consist of the vertices of a certain tessellation with prisms and possibly other points. Roughly speaking, the parameter $a$ determines the size of the basis of those prisms and the parameter $b$ determine the height of each prism.

There are some obvious restrictions on the parameters $a$ and $b$. If we want the vertex set of $\PL$ to be finite, the parameter $a$ must be rational. It is not hard to see that for planar polyhedra there are no restrictions for $b$. The same technique used in pure polyhedra applies to determine the possible values for $a$ for every planar polyhedron. Those values are listed in Table\nobreakspace \ref {tab:parametersPlanar}. It is necessary forbid small values since they will break the diamond condition in $\PL$. In this table we list only one polyhedron by Petrie pair since the same parameters apply for Petrials of planar polyhedra.

\begin{table}
\begin{center}
{%
\begin{tabular}{|cll|}
\hline
Polyhedron & \multicolumn{2}{c|}{Possible values of $a$}\\
\hline
$\{4,4\}$ & $a \in \bN \sm \{1\}$  & if $\bfL D^{-1} = \Cl$\\
 & $a \in \frac{1}{2}(\bN \sm \{1,2\})$ & if $\bfL D^{-1} = \Fcl$\\
 & $a \in \bN \sm \{1\}$& if $\bfL D^{-1} = \Bcl$\\
 & $a \in \frac{1}{2}(\bN \sm \{1,2\})$ & if $\bfL D^{-1} = \langle \Scl, t_{3} \rangle$\\ \hline
$\{3,6\}$ & $a \in \bN \sm \{1\}$ & if $\bfL D^{-1} = \langle \Tl, t_{3} \rangle$  \\
 &  $a \in \frac{1}{3}(\bN \sm \{1\})$ & if $\bfL D^{-1} = \langle \Tcl, t_{3} \rangle$\\ \hline
$\{6,3\}$ & $a \in \frac{1}{2}(\bN \sm \{1,2\})$ & if $\bfL D^{-1} = \langle \Tl, t_{3} \rangle$\\
 & $a \in \frac{1}{2}\bN$ & if $\bfL D^{-1} = \langle \Tcl, t_{3} \rangle$\\
\hline
\end{tabular}
}
\caption{Parameters for planar apeirohedra} \label{tab:parametersPlanar}
\end{center}
\end{table}

Before giving the parameter for blended polyhedra we discuss some details about them. Recall that realizations of blended regular polyhedra depend on a parameter $\alpha$ (see Section\nobreakspace \ref {sec:basics}). In order to make calculations simpler, in this subsection this parameter will be assumed to be $1$ for $\{4,4\} \# \{\ \}$, $\{4,4\} \# \{\infty\}$, $\{3,6\} \# \{\ \}$ and $\{3,6\} \# \{\infty\}$ and $\frac{1}{3}$ for $\{6,3\}\# \{\ \}$ and $\{6,3\}\# \{\infty \}$. Calculations made in this subsection might be easily modified to arbitrary $\alpha$.

Following the notation of Theorem\nobreakspace \ref {thm:infinite} we have to determine the pairs $(a,b)$ of parameters that induce toroidal realizations. As mentioned before, there are three kinds of undesired identifications, those that identify a vertex with an interior point of an edge, those that identify a vertex with a non-translate vertex making the vertex figures no to be cycles, and those that identify two non-translate vertices of the same face. 

The first consideration is the same mentioned before for planar polyhedra: the parameter $a$ must be rational, otherwise the vertex set of $\PL$ is not finite.

Observe that if $\P$ is a blended polyhedra of type $\Q \# \{\ \}$ and $b > 1$ then any choice of $a \in \Q$ will work, therefore we only have to consider $0< b \leq 1$. Since $a$ is rational and the slope of any edge is also rational, if $b \in (0.1] \sm \bQ$ then the pair $(a,b)$ will not produce undesired identifications. If $\P$ is of type $\Q \# \{\infty\}$, $b$ must be also rational. 

Given the considerations mentioned before, it only remains to determine which of the pairs of parameters $(a,b)= \left(\frac{p}{q}, \frac{r}{s}\right)$, with $p,q,r,s \in \bZ$ and $\gcd(p,q)=\gcd(r,s) = 1$, induce toroidal realizations (bear in mind that if $\P$ is of type $\Q \# \{\ \}$ we may assume $r \leq s$). In Tables\nobreakspace  \ref {tab:parametersBlended1} to\nobreakspace  \ref {tab:parametersBlended4}  we give necessary and sufficient conditions over $p$, $q$, $r$, $s$ and $g = \gcd(q,s)$ for every lattice and every blended polyhedron (up to Petri operation) to produce identifications between a vertex or an interior point of an edge or between two non-translate vertices of $\P$. The calculations to determine the values are long but straightforward and usually involve determining the existence of a certain solution of a system of linear diophantine equations.

We still have to consider those identifications that occur between pairs of translate vertices but identify vertices in the same face, since this situation might break diamond condition in $\PL$. This phenomenon is analogous to that in the map $\{4,4\}_{(1,1)}$ (see \cite[Section 1D]{ARP}). Observe that in the situation just described diamond condition holds if and only if the translation that identifies two vertices in the same face maps the face to itself. In particular, if the faces are finite, then no two vertices in the same face can be identified.

If $\P = \{4,4\} \# \{\ \}$ or its Petrial the only case that has to be considered is whenever the vector $(1,1,0)$ belongs to $\L$, i.e. $p=1$ if $\L D^{-1} \in \{\cl, \fcl,  \scl \+ (0,0,1)\}$ and $p =1$, $2 \mid q$ if $\L D^{-1} = \bcl$. As mentioned before, if $\P = \{4,4\} \# \{\ \}$ this breaks diamond condition. If $\P = \{\infty,4\}_{4} \# \{\ \}$ the resulting polyhedron $\PL$ is the unique polyhedron with two vertices of type $\{2,4\}$.

If $\P = \{4,4\} \# \{\infty \}$ we shall forbid $p=1$ if $\L D^{-1} \in \{\cl, \fcl,  \scl \+ (0,0,1)\}$ and $p =1$, $2 \mid q$ if $\L D^{-1} = \bcl$ since this implies that $(1,1,2r) \in \L$, which identifies a vertex with a non-translate vertex of the same face. In $\{\infty,4\}_{4} \# \{\ \}$ there are no extra forbidden values.

For $\{3,6\} \# \{\infty \}$  we must forbid the occurrence of $p=1$ and $3 \nmid r$ simultaneously if $\L D^{-1} = \tl \+ (0,0,1)$ and $p=1$, $3 \nmid r$, $3 \mid q$ whenever $\L$ is any other of the lattices preserved by $G_{o}(\{3,6\} \# \{\infty \})$ since this will imply that $(1,0,3k+1) \in \L$ for some $k \in \bZ$ and this vector produces an identification that breaks diamond condition in $\PL$. For $\{\infty,6\}_{3} \# \{\infty \}$ we must forbid $p=1$ and $2 \nmid r$ since this will imply that $(1, \frac{3r-1}{2},3r) \in \L$ and the corresponding translation identifies two vertices in the same face but do not preserves the face. 

If $\{6,3\} \# \{\ \}$ and $\L D^{-1} \in \{\tl \+ e_3, 3 \tl \+ (u_1 + u_2+e_3), 3 \tl \+ (u_1 + u_2-e_3)\}$, with $\{u_1,u_2\}$ the basis of $\tl$ given at the end of Section\nobreakspace \ref {sec:lattices}, we must forbid $p=1$. If $\L D^{-1} = \tcl \+ e_3$ there is no need to forbid anything else. It is not necessary to forbid any other set of parameters for the polyhedron $\{\infty,3\}_{6} \# \{\ \}$.

If $\L D^{-1} = \tl \+ e_3$ and $\P = \{3,6\} \# \{\ \}$ the only extra condition necessary to avoid undesired identifications is $p=1$; if $\L D^{-1} = \tcl \+ e_3$ there is no need to add extra conditions. No additional condition is needed if $\P = \{\infty,6\}_{3} \# \{\ \}$. 

With the polyhedron $\{6,3\} \# \{\infty\}$ the only extra condition of $p \neq 1$ must be added whenever $\L D^{-1} = \tl + e_3$, since $p=1$ implies that $(1,0,2r) \in \L$ and the corresponding translation identifies two vertices in the same face and do not preserve the face. No additional condition is needed when $\L D^{-1} = \tcl + e_3$. There is no need to add any condition if $P = \{ \infty,3\}_{6} \# \{\infty\}$.

\begin{table}
\begin{center}
 \resizebox{\textwidth}{!} {
{
\newcommand{\mc}[3]{\multicolumn{#1}{#2}{#3}}
\begin{tabular}{|c|p{.2\textwidth}|p{.2\textwidth}|p{.2\textwidth}|p{.2\textwidth}|}
	\hline
	$\P$ &  \mc{1}{c|}{$\Cl$}& \mc{1}{c|}{$\Fcl$} &\mc{1}{c|}{$\Bcl$} &\mc{1}{c|}{$\langle \Scl, t_3 \rangle$} \\ \hline
	$\{4,4\} \# \{\ \}$ &
	\begin{minipage}[t]{\linewidth}
		\begin{itemize}[leftmargin=10pt, noitemsep, topsep=0pt]%
			\item $2 \nmid p$, $2 \nmid r$ and $r \leq g$.
			\item $2 \nmid p$, $2 \mid r$ and $r < g$.
			\item $2 \mid p$, $2 \nmid r$ and $2r < g$.
			\item $2 \mid p$, $2 \mid r$ and $r < g$.
		\end{itemize}
	\end{minipage} &
	\begin{minipage}[t]{\linewidth}
		\begin{itemize}[leftmargin=10pt, noitemsep, topsep=0pt]
			\item $2 \nmid p$, $2 \nmid r$, $2 \mid \frac{q}{g}$ and $2r \leq g$.
			\item $2 \nmid p$, $2 \nmid r$, $2 \nmid \frac{q}{g}$, $2 \nmid \frac{s}{g}$  and $r \leq g$.
			\item $2 \nmid p$, $2 \nmid r$, $2 \nmid \frac{q}{g}$, $2 \mid \frac{s}{g}$, $2 \nmid g$  and $r \leq g$.
			\item $2 \nmid p$, $2 \nmid r$, $2 \nmid \frac{q}{g}$, $2 \mid \frac{s}{g}$, $2 \mid g$  and $2r \leq g$.
			\item $2 \nmid p$, $2 \mid r$ and $2r < g$.
			\item $2 \mid p$, $2 \nmid r$ and $2r < g$.
			\item $2 \mid p$, $2 \mid r$ and $r < g$.
		\end{itemize}
	\end{minipage}&
	\begin{minipage}[t]{\linewidth}
		\begin{itemize}[leftmargin=10pt, noitemsep]
			\item $2 \nmid p$, $2 \nmid r$ and $2r \leq g$.
			\item $2 \nmid p$, $2 \mid r$, $2 \mid q$ and $2r < g$.
			\item $2 \nmid p$, $2 \mid r$, $2 \nmid q$ and $r < g$.
			\item $2 \mid p$, $2 \nmid r$ and $2r < g$.
			\item $2 \mid p$, $2 \mid r$ and $2r < g$.
		\end{itemize}
	\end{minipage}&
	\begin{minipage}[t]{\linewidth}
		\begin{itemize}[leftmargin=10pt, noitemsep]
			\item $2 \nmid p$, $2 \nmid r$, $2 \nmid \frac{q}{g}$ and $2r \leq g$.
			\item $2 \nmid p$, $2 \nmid r$, $2 \mid \frac{q}{g}$ and $r \leq g$.
			\item $2 \nmid p$, $2 \mid r$ and $r < g$.
			\item $2 \mid p$, $2 \nmid r$ and $2r < g$.
			\item $2 \mid p$, $2 \mid r$ and $r < g$.
		\end{itemize}
	\end{minipage} \\ \hline

	$\{4,4\}\# \{\infty\}$ &%
	\begin{minipage}[t]{\linewidth}
		\begin{itemize}[leftmargin=10pt, noitemsep]
			\item $g=1$, $2 \nmid p$ and $2 \nmid r$.
			\item $g=3$ and do not occur that $r \equiv \pm 4 \pmod{12}$ and $2\mid p$.
			\item $g \notin \{1,3\}$.
		\end{itemize}
	\end{minipage}&
	\begin{minipage}[t]{\linewidth}
		\begin{itemize}[leftmargin=10pt, noitemsep]
			\item $g=1$, $2 \nmid p$, $2 \nmid r$, $2 \nmid q$, $2 \nmid s$.
			\item $g=3$ and none of the following occur: \begin{itemize}[leftmargin=12pt, noitemsep] \item $r \equiv \pm 4 \pmod{12}$ and $2\mid p$. \item $r \equiv \pm 2 \pmod{12}$ and $2 \nmid pq$. \end{itemize}
			\item $g \notin \{1,3\}$.
		\end{itemize}
	\end{minipage}&
	\begin{minipage}[t]{\linewidth}
		\begin{itemize}[leftmargin=10pt, noitemsep]
			\item $g=3$ and do not occur that $r \equiv \pm 4 \pmod{12}$ and $2\mid p$.
			\item $g \notin \{1,3\}$.
		\end{itemize}
	\end{minipage}&
	\begin{minipage}[t]{\linewidth}
		\begin{itemize}[leftmargin=10pt, noitemsep]
			\item $g=1$, $2 \nmid p$, $2 \nmid r$ and $2 \mid q$.
			\item $g=3$ and do not occur that $r \equiv \pm 4 \pmod{12}$ and $2\mid p$.
			\item $g \notin \{1,3\}$.
		\end{itemize}
	\end{minipage} \\ \hline
\end{tabular}

}

 }
\caption{Parameters for $\{4,4\} \# \{\ \}$, $\{4,4\} \# \{\infty\}$ and their Petrials} \label{tab:parametersBlended1}
\end{center}
\end{table}

\begin{table}
\begin{center}
 \resizebox{\textwidth}{!} {
{
\newcommand{\mc}[3]{\multicolumn{#1}{#2}{#3}}
\begin{tabular}{|c|p{.2\textwidth}|p{.2\textwidth}|p{.2\textwidth}|p{.2\textwidth}|}
	\hline
	$\P$ &  \mc{1}{c|}{$\left\langle \Tl, t_3 \right\rangle$}& \mc{1}{c|}{$\left\langle \Tcl, t_3 \right\rangle $} &\mc{1}{c|}{$\left\langle \Tcl, t_1 \right\rangle$} &\mc{1}{c|}{$\left\langle \Tcl, t_{2}\right\rangle$} \\ \hline
	$\{3,6\} \# \{\infty\}$ &
	\begin{minipage}[t]{\linewidth}
		\begin{itemize}[leftmargin=10pt, noitemsep, topsep=0pt]%
			\item $g \neq 1$.
		\end{itemize}
	\end{minipage}&	
	\begin{minipage}[t]{\linewidth}
		\begin{itemize}[leftmargin=10pt, noitemsep, topsep=0pt]
			\item $g \not\in \{1,3\}$.
		\end{itemize}
	\end{minipage}&	
	\begin{minipage}[t]{\linewidth}
		\begin{itemize}[leftmargin=10pt, noitemsep]
			\item $g=3$ and $\frac{ps}{3} \equiv  \frac{rq}{3} \pmod{3} $
			\item $g \not\in \{1,3\}$.
		\end{itemize}
	\end{minipage}&	
	\begin{minipage}[t]{\linewidth}
		\begin{itemize}[leftmargin=10pt, noitemsep]
			\item $g=3$ and $\frac{ps}{3} \equiv  -\frac{rq}{3} \pmod{3} $
			\item $g \not\in \{1,3\}$.
		\end{itemize}
	\end{minipage} \\ \hline
\end{tabular}

}
 }
\caption{Parameters for $\{3,6\} \# \{\infty\}$ and its Petrial.} \label{tab:parametersBlended2}
\end{center}
\end{table}

\begin{table}
\begin{center}
 \resizebox{\textwidth}{!} {

 {
\newcommand{\mc}[3]{\multicolumn{#1}{#2}{#3}}
\begin{tabular}{|c|p{.2\textwidth}|p{.2\textwidth}|p{.2\textwidth}|p{.2\textwidth}|}
	\hline 
	$\P$ &  \mc{1}{c|}{$\left\langle \Tl, t_3 \right\rangle$}& \mc{1}{c|}{$\left\langle \Tcl, t_3 \right\rangle $} &\mc{1}{c|}{$\left\langle 3\Tl, t_1 \right\rangle$} &\mc{1}{c|}{$\left\langle 3\Tcl, t_{2}\right\rangle$} \\ \hline
	$\{6,3\} \# \{\ \}$ &
	\begin{minipage}[t]{\linewidth}
		\begin{itemize}[leftmargin=10pt, noitemsep, topsep=0pt]%
			\item $3r \leq g$.
		\end{itemize}
	\end{minipage}&	
	\begin{minipage}[t]{\linewidth}
		\begin{itemize}[leftmargin=10pt, noitemsep, topsep=0pt]
			\item $3r \leq g$.
		\end{itemize}
	\end{minipage}&	
	\begin{minipage}[t]{\linewidth}
		\begin{itemize}[leftmargin=10pt, noitemsep]
			\item $3 \nmid q$ and $3r < g$.
			\item $3 \mid \frac{q}{g}$ and $9r \leq g$.
			\item $3 \mid q $, $3 \nmid \frac{q}{g}$, $3 \mid \frac{s}{g}$ and $9r \leq g$.
			\item $3 \mid q $, $3 \nmid \frac{q}{g}$, $3 \nmid \frac{s}{g}$, $\frac{rq}{g} \equiv \frac{sp}{g}$ and $3r \leq g$.
			\item $3 \mid q $, $3 \nmid \frac{q}{g}$, $3 \nmid \frac{s}{g}$, $\frac{rq}{g} \not\equiv \frac{sp}{g}$ and $9r \leq g$.
		\end{itemize}
	\end{minipage}&	
	\begin{minipage}[t]{\linewidth}
		\begin{itemize}[leftmargin=10pt, noitemsep]
			\item $3 \nmid q$ and $3r < g$.
			\item $3 \mid \frac{q}{g}$ and $9r \leq g$.
			\item $3 \mid q $, $3 \nmid \frac{q}{g}$, $3 \mid \frac{s}{g}$ and $9r \leq g$.
			\item $3 \mid q $, $3 \nmid \frac{q}{g}$, $3 \nmid \frac{s}{g}$, $\frac{rq}{g} \equiv -\frac{sp}{g}$ and $3r \leq g$.
			\item $3 \mid q $, $3 \nmid \frac{q}{g}$, $3 \nmid \frac{s}{g}$, $\frac{rq}{g} \not\equiv -\frac{sp}{g}$ and $9r \leq g$.
		\end{itemize}
	\end{minipage} \\ \hline
\end{tabular}

}
 
 }
\caption{Parameters for $\{6,3\} \# \{\ \}$ and its Petrial.} \label{tab:parametersBlended3}
\end{center}
\end{table}

\begin{table}
\begin{center}
{
\newcommand{\mc}[3]{\multicolumn{#1}{#2}{#3}}
\begin{tabular}{|c|p{.35\textwidth}|p{.35\textwidth}|}
	\hline 
	$\P$ &  \mc{1}{c|}{$\left\langle \Tl, t_3 \right\rangle$}& \mc{1}{c|}{$\left\langle \Tcl, t_3 \right\rangle $}  \\ \hline
	$\{3,6\} \# \{\ \}$ &
	\begin{minipage}[t]{\linewidth}
		\begin{itemize}[leftmargin=10pt, noitemsep, topsep=0pt]%
			\item $r < g$.
		\end{itemize}
	\end{minipage}&	
	\begin{minipage}[t]{\linewidth}
		\begin{itemize}[leftmargin=10pt, noitemsep, topsep=0pt]
			\item $3 \nmid q$ and $r < g$.
			\item $3 \mid \frac{q}{g}$ and $r < g$.
			\item $3 \nmid \frac{q}{g}$, $3 \mid g$ and $3r < g$.
		\end{itemize} 
	\end{minipage} \\ \hline
	$\{6,3\} \# \{\infty \}$ &
	\begin{minipage}[t]{\linewidth}
		\begin{itemize}[leftmargin=10pt, noitemsep, topsep=0pt]%
			\item $g=1$, $3 \nmid r$ and $2 \nmid r$.
			\item $g=2$, $3 \nmid p$ and $4 \nmid r$.
			\item $g=5$ and $3 \nmid p$.
			\item $g=5$ and $2 \nmid r$.
			\item $g \notin \{1,2,5\}$.
		\end{itemize}
	\end{minipage}&	
	\begin{minipage}[t]{\linewidth}
		\begin{itemize}[leftmargin=10pt, noitemsep, topsep=0pt]
			\item $g=1$, $3 \nmid r$ and $2 \nmid r$.
			\item $g=2$, $3 \nmid p$ and $4 \nmid r$.
			\item $g=5$ and $3 \nmid p$.
			\item $g=5$ and $2 \nmid r$.
			\item $g \notin \{1,2,5\}$.
		\end{itemize} 
	\end{minipage} \\ \hline
\end{tabular}

}

\caption{Parameters for $\{3,6\} \# \{\ \}$, $\{6,3\} \# \{\infty \}$ and their Petrials.} \label{tab:parametersBlended4}
\end{center}
\end{table}

\FloatBarrier

\subsection{Completeness of the list}

The aim of this subsection is to prove that the regular polyhedra in $\TT$ described above are all regular polyhedra in $\TT$. Roughly speaking, given a regular polyhedron $\PL$ in $\TTL$ we will construct a regular polyhedron $\P$ in $\E$ such that the quotient map $\piL$ maps $\P$ to $\PL$. We will not give full detail of the proofs here but most of them are consequence the theory developed in \cite[Chapter 8]{ratcliffeHypManifolds}, in particular of Theorem\nobreakspace \ref {thm:isomTorus} and the following lemma.

\begin{lem} \label{lem:localIsom}
	If $\bfL$ is a lattice group, the quotient map $\piL$ is a local isometry.
\end{lem}

Now we briefly explain how to construct the regular polyhedron $\P$ in $\E$ given a regular polyhedron $\PL$ in $\TTL$. Let $\{F_0, F_1, F_2\}$ be the base flag of $\PL$ with $rk(F_{i}) = i$, $i \in \{0,1,2\}$. Without loss of generality we may assume that $F_0 = o \bfL$, equivalently, that $o$ projects to $F_{0}$. Let $V$ be the vertex set of $\PL$, and consider $W_{0}= \pi^{-1}_{\bfL} (V) = \bigcup_{x \in V} \{\pi^{-1}_{\bfL} (x)\}$. Let $p$ be a point in $\E$ such that $p\bfL$ is the other point of $F_1$ different from $o$. By arc-lifting theorem on covering maps there exists a unique curve $\gamma: [0,1] \to \E$ that projects to $F_1$ and such that $\gamma(0) = o $ and $\gamma(1) \in p \bfL$; without loss of generality we may assume that $\gamma(1)=p$. By Lemma\nobreakspace \ref {lem:localIsom}, $\gamma$ must be the line segment joining $o$ and $p$. Construct a family $W_{1}$ of line segments in $\E$ with endpoints in $W_{0}$ as follows: take any edge $G_{1}$ of $\PL$ and an isometry $\bar{S}$ of $\TTL$ mapping $F_{1}$ to $G_{1}$. By Theorem\nobreakspace \ref {thm:isomTorus} there is a family of isometries of $\E$ that project to $\bar{S}$ and any two of them differ by an element of $\bfL$, consider the image of $\gamma$ under all those isometries; $W_{1}$ consists of all line segments constructed this way varying $G_{1}$ over the edges of $\PL$. Observe that these line segments are precisely the possible liftings of the edges of $\PL$, and that $\gamma S$ projects to $F_1 \bar{S}$ if and only if $S$ projects to $\bar{S}$.

The points in $W_{0}$ together with the line segments in $W_{1}$ give a (geometric) graph in $\E$. If it is connected then this graph will be the $1$-skeleton of $\P$, if not, just take as $\sk(\P)$ the connected component of this graph that contains the origin $o$. 

The faces of $\P$ are induced by $\PL$ in the following way: consider the face $F_2$, it may be thought as a subgraph of $\sk(\PL)$ and hence, as a loop in $\TTL$. If this loop is contractible, then there is a unique cycle $\omega$ of $\sk(\P)$ containing $o$ that projects to $F_2$. If $F_2$ is not contractible there is a unique path $P$ in $\sk(\P)$ that starts in $o$, ends in $ot$ for some $t \in \bfL$ and projects to $F_2$. The orbit of $P$ under $\langle t \rangle$ gives place to an helix or a zigzag $\omega$ in $sk(\P)$. In any case, $\omega$ will define a face of $\P$. The other faces of $\P$ will be constructed from $\omega$ in the analogous way as the other edges of $\P$ were constructed from $\gamma$. Again $\omega S$ projects to $F_2 \bar{S}$ if and only if $S$ projects to $\bar{S}$.

Now we have constructed a connected graph $\sk(\P)$ and a family of subgraphs of $\sk(\P)$ such that every subgraph is either an helix or a cycle. It is clear that every edge belong to at least one face of $\P$ and since $\pi_{\bfL}$ is a local isometry, the vertex figure of a vertex $x$ of $\P$ must be a cycle since the vertex figure of $\pi_{\bfL}(x)$ was a cycle. Therefore the structure $\P$ induced by this construction must be a polyhedron. By construction it is clear that $\pi_{\bfL}$ induces a graph epimorphism from $\sk(\P)$ to $\sk(\PL)$. It only rest to prove that $\P$ is regular.

Assume that $\bar{R}_{0}$, $ \bar{R}_{0} $ and $\bar{R}_{2}$ act as distinguished generators in $G_{\bfL}(\PL)$ for the base flag $\{F_0,F_1,F_2\}$. By Theorem\nobreakspace \ref {thm:isomTorus} there exist unique isometries $R_0$, $R_1$, $R_2$ such that $R_i$ projects to $\bar{R}_{i}$, $i \in \{0,1,2\}$, $R_1$ and $R_2$ fix the origin $o$ and $R_{0}$ fix the midpoint of $\gamma$. Observe that $R_0$ maps $o$ to a translate of $p$, the other vertex in $\gamma$, but since it maps $\gamma$ to itself, it must map $o$ to $p$. Similarly, $R_0$ maps $\omega$ to a translate of itself but since it fixes $\gamma$ it must map $\omega$ to itself. Since $\bar{R}_{1}$ fixes $F_2$, $R_1$ maps $\omega$ to a translate of itself but since $R_{1}$ fixes $o$ then it fixes $\omega$; $R_1$ must map $\gamma$ to the unique edge of $\P$ containing $o$ that projects to $F_1 \bar{R}_{1}$. Similarly, $R_2$ must fix $\gamma$ and map $\omega$ to the unique face of $\P$ containing $o$ that projects to $F_2 \bar{R}_{2}$. Therefore, $R_0$, $R_1$ and $R_2$ act as distinguished generators of $G(\P)$ for the base flag $\{o,\gamma, \omega \}$, thus $\P$ is regular. We have proved the following result.

\begin{thm}\label{thm:lifting}
	Let $\bfL$ be a lattice group and $\PL$ a regular polyhedron in $\TTL$. There exists a regular polyhedron $\P$ in $\E$ such that $\pi_{\bfL}(\P) = \PL$.
\end{thm}

Theorems\nobreakspace  \ref {thm:Ico} to\nobreakspace  \ref {thm:infinite}  and\nobreakspace  \ref {thm:lifting} complete the classification of regular polyhedra in the $3$-torus.

\bibliographystyle{plain}
\bibliography{bib}

\end{document}